\documentclass[journal ]{new-aiaa}
\usepackage[utf8]{inputenc}
\usepackage{textcomp}

\usepackage{graphics} % for pdf, bitmapped graphics files
\usepackage{epsfig} % for postscript graphics files
\usepackage{amsmath}
\usepackage[version=4]{mhchem}
\usepackage{siunitx}
\usepackage{enumitem}
\usepackage{longtable,tabularx}
\usepackage{multirow}
\usepackage{subcaption}
\makeatletter
\newcommand{\thickhline}{%
    \noalign {\ifnum 0=`}\fi \hrule height 1pt
    \futurelet \reserved@a \@xhline
}

\title{Analysis of Cooperative and Non-Cooperative Architectures for Multi-Plane On-Orbit Refueling}

\author{Kosuke Ikeya \footnote{Ph.D. Student, Department of Earth Science and Engineering. Student Member AIAA.}}
\affil{Imperial College London, London, SW7 2AZ, United Kingdom}
\author{Koki Ho\footnote{Associate Professor, Daniel Guggenheim School of Aerospace Engineering. Senior Member AIAA.}}
\affil{Georgia Institute of Technology, Atlanta, Georgia, 30332}

\begin{document}

\maketitle

\begin{abstract}
% As many satellite constellations are proposed, deployed, and operated, their maintenance becomes increasingly important to provide satisfactory services; therefore, on-orbit refueling to spacecraft has become one of the most promising technologies that realize more sustainable satellite constellations. This paper develops an analytical model to examine two types of mission architectures for multi-target on-orbit refueling missions: a non-cooperative architecture and a cooperative architecture. While in the (rather conventional) non-cooperative refueling architecture, a servicer spacecraft visits passive targets one by one, in the cooperative refueling architecture, both the servicer and the targets can actively maneuver to complete refueling cooperatively.
% This paper analytically compares the fuel mass required in each architecture to support the decision-making process of mission architects. Notably, a crossover point of a mass ratio between the servicer and the target where the suggested mission architecture changes from the non-cooperative to the cooperative is analytically derived.
% The effects of mission parameters such as the number of targets and inclination angle on this ratio are also examined through a case study of multi-plane multi-target on-orbit refueling in low Earth orbits.
% The result suggests that cooperative architectures require less fuel when the actual mass ratio between the servicer and the target is larger than the critical mass ratio. {\color{red}The sensitivity of this ratio against key mission parameters is analyzed. }

As many satellite constellations are proposed, deployed, and operated, their maintenance becomes increasingly important to provide satisfactory services; therefore, on-orbit refueling to spacecraft has become one of the most promising technologies for realizing more sustainable space development. This paper develops an analytical model to examine two types of mission architectures for multi-target on-orbit refueling missions: a non-cooperative architecture and a cooperative architecture. In the (rather conventional) non-cooperative refueling architecture, a servicer spacecraft visits passive targets one by one, whereas, in the cooperative refueling architecture, both the servicer and the targets can actively maneuver to complete refueling cooperatively.
This paper analytically compares the fuel mass required in each architecture to support the decision-making process of mission architects. Furthermore, the condition under which the cooperative architecture becomes more efficient than the non-cooperative architecture is analytically derived. 
The sensitivities of this condition against key mission parameters, such as the number of targets and their inclination, are also analyzed through a case study of multi-plane multi-target on-orbit refueling in low Earth orbits.

\end{abstract}

\section*{Nomenclature}

% \noindent(Nomenclature entries should have the units identified)

{\renewcommand\arraystretch{1.0}
\noindent\begin{longtable*}{@{}l @{\quad=\quad} l@{}}
$g_0$ & standard acceleration due to gravity\\
$I_\mathrm{sp}$ & specific impulse\\
$i$ & inclination\\
$m$ & general mass\\
% $m_\mathrm{s}$ & servicer mass\\
% $m_\mathrm{s,I}$ & servicer initial mass at the beginning of the whole sequence\\
% $m_\mathrm{s,F}$ & servicer final mass at the end of the whole sequence\\
% $m_\mathrm{t,T}$ & target initial mass at the beginning of the whole sequence\\
$m_\mathrm{r}$ & total refuel mass\\
$m_\mathrm{req}$ & required refuel mass\\
% $m_\mathrm{t,i}$ & target mass at inbound transfer\\
% $m_\mathrm{t,o}$ & target mass at outbound transfer\\
$n$  & number of targets \\
$u$ & relative argument of latitude\\
$\alpha$ & critical mass ratio between a servicer and targets where the suggested mission architecture changes\\
$\Delta v$ & change in velocity\\
% $a$ &    cylinder diameter \\
% $C_p$& pressure coefficient \\
% $Cx$ & force coefficient in the \textit{x} direction \\
% $Cy$ & force coefficient in the \textit{y} direction \\
% c   & chord \\
% d$t$ & time step \\
% $Fx$ & $X$ component of the resultant pressure force acting on the vehicle \\
% $Fy$ & $Y$ component of the resultant pressure force acting on the vehicle \\
% $f, g$   & generic functions \\
% $h$  & height \\
% $i$  & time index during navigation \\
% $j$  & waypoint index \\
% $K$  & trailing-edge (TE) nondimensional angular deflection rate\\
% $\Theta$ & boundary-layer momentum thickness\\
% $\rho$ & density\\
\multicolumn{2}{@{}l}{Subscripts}\\
c & cooperative case\\
% d & dry mass\\
F & mass after a whole campaign\\
f & mass after an orbital transfer\\
I & mass before a whole campaign\\
in & mass at inbound transfer\\
n & non-cooperative case\\
out & mass at outbound transfer\\
s & servicer\\
t & target\\
0 & mass before an orbital transfer
% $G$ & generator body\\
% iso	& waypoint index
\end{longtable*}}

\section{Introduction}
\lettrine{S}{atellite} constellations have been proposed, researched, and operated for various objectives such as Earth observation, navigation, and communication \cite{gree1989, fore2017, bosh2014, lee2020, ande2022}.
For instance, Planet operates a constellation, PlanetScope, to obtain images of Earth's land frequently. These images acquired from satellite constellations are, for example, used for remote sensing \cite{roy2021}.
The Global Positioning System (GPS), which is now broadly used to navigate many devices such as cars and smartphones, also utilizes a satellite constellation of about 30 satellites in six different orbital planes in medium Earth orbit.
Among communication constellations, constellations for internet access, such as Starlink operated by SpaceX, are becoming popular. SpaceX plans to deploy about 12,000 satellites in total in Low Earth Orbit (LEO) to provide low-latency internet access in many countries. These services from constellations have become essential to our daily lives and emergencies, such as disaster operations management and Internet access during a war.

To keep providing satisfactory performances to users, maintaining these constellations is crucial.
% , and several approaches have been proposed.
Replacing malfunctioning satellites has been proposed and researched as one maintenance  approach \cite{corn1999, jako2019}. 
Cornara et al.~\cite{corn1999} analyzed three different replacement strategies: launching new satellites on-demand, placing spare satellites in close orbits, and placing these spare satellites in parking orbits. Jakob et al.~\cite{jako2019} proposed an analytical model of the parking orbits strategy by leveraging a multi-Echelon inventory policy.

Another potential maintenance approach is to service these satellites on orbit (On-Orbit Servicing, OOS)\cite{luu2022, hall1999}.
The validity of the general OOS concept has been shown in Refs.~\cite{long2007,hatt2022, luu2022}. For instance, Hatty~\cite{hatt2022} analyzed the economical viability of extending the operational period of Landsat 7 by a robotic spacecraft, OSAM-1, and compared it to the cost of replacing the satellite. Luu and Hastings~\cite{luu2022} focused on LEO constellations and explored a tradespace where OOS could be beneficial.
Some flight demonstrations for essential technologies for OOS were also performed \cite{inab2000,ogil2008}. ETS-VII developed and operated by NASDA (now merged into JAXA) demonstrated a berthing operation to capture a passive free-flying target {\cite{oda1999, inab2000}}. The Orbital Express Demonstration System mission managed by DARPA also demonstrated capturing and berthing of a free-flying target followed by transferring of Orbital Replacement Units {\cite{ogil2008}}.
Furthermore, successful dockings of Mission Extension Vehicles to actual commercial satellites, Intelsats, demonstrates the feasibility of prolonging the operational period of satellites\cite{nort2020, nort2021}, As these examples show, OOS is becoming even more promising.

% \iffalsen Recent studies in OOS have shown the validity of the concept  and proposed many types of missions that can exploit OOS,  such as  active debris removal and\cite{hong2018}
% life extension of satellites\fi

An OOS mission to a satellite constellation can be considered a multi-target OOS mission. 
Analyses and optimizations of various kinds of multi-target OOS missions have been extensively researched \cite{jean2006, meng2019, hong2018,sart2017,sart2021, sart2022, ho2020}.
Bourjolly et al.~\cite{jean2006} minimized the total energy or time required for an OOS mission by solving a time-dependent, moving-target traveling salesman problem.
Meng et al.~\cite{meng2019} proposed a new system architecture for on-orbit refueling, “1+N”, consisting of one fuel storage station and $N$ refueling vehicles, and a multi-objective optimization was conducted to find an optimal deployment strategy.
Shen et al.~\cite{hong2018} also conducted a multi-objective optimization, but for a multiple debris removal mission.
Sarton du Jonchay and Ho \cite{sart2017} quantitatively compared the responsiveness of two OOS infrastructures, the one "With Depot" and the one "Without Depot." 
Ho et al.~\cite{ho2020} further developed a semi-analytical OOS model with a depot considering the queue of services and the capacity of the depot.
In Ref.~\cite{sart2021}, Sarton du Jonchay et al. developed an optimization framework to help decision-makers to make operational and strategic plans under uncertainties in service demand.
Sarton du Jonchay et al.~\cite{sart2022} further extended this framework to handle OOS with low-thrust orbital transfers.

These studies reviewed above have employed mission architectures with one active servicer or multiple active servicers and passive targets (i.e., non-cooperative OOS architectures). Another approach, peer-to-peer (P2P) OOS (especially on-orbit refueling) strategies have also been researched \cite{shen2005, sala2006, dutt2008, dutt2010, dutt2012}. In the later egalitarian P2P refueling studies \cite{dutt2010, dutt2012}, all spacecraft in a constellation are assumed to be capable of conducting orbital transfers and refueling another spacecraft.
Dutta and Tsiotras~\cite{dutt2010} minimized the fuel consumption during the orbital transfers proposing a network flow formulation of this problem. Dutta et al.~\cite{dutt2012} further extended this strategy to low-thrust maneuvers, developed a solver to solve the optimization problem, and demonstrated it through two mission scenarios.

In addition to these two architectures in on-orbit constellation servicing, non-cooperative OOS and P2P OOS, a new architecture, cooperative OOS, has been recently proposed \cite{du2015,zhao2017, cox2022}. This architecture is a mixture of non-cooperative and P2P OOS allowing maneuvering by any spacecraft in the architecture while targets (spacecraft to be serviced) cannot service any other spacecraft.
Zhao et al.~\cite{zhao2017} analytically derived an expression of the total propellant mass for cooperative architectures and compared it to the cost with a traditional non-cooperative architecture considering the $J_2$ perturbation. They showed the benefit of employing this architecture through a case study of a coplanar on-orbit refueling mission.  
% {\color{red}Their research scope is limited to coplanar OOS missions}.
Cox et al.~\cite{cox2022} employed the consensus-based bundle algorithm (CBBA) for the initial planning and scheduling of the cooperative OOS.
Refs.~\cite{du2015,zhao2017} compared the cooperative and the conventional non-cooperative architectures and showed the potential less fuel required in the cooperative architectures.
However, the effect of the target constellation parameters such as orbital elements and the number of targets on the required fuel mass in both architectures is yet to be discovered.
To provide a more general and comprehensive comparison between cooperative and non-cooperative architectures to help mission architects, this paper analyses non-coplanar multi-target OOS missions with multiple mission architectures, and addresses the conditions that make one architecture better than the other.

The rest of the paper is organized as follows: mathematical derivation of the initial mass of a servicer in a cooperative architecture is introduced in Sec.~\ref{sec:problem_formulation}. Based on the derived initial mass, Sec.~\ref{sec:analytical_comparison} compares cooperative and non-cooperative architecture analytically. Sec.~\ref{sec:case} discusses the value of cooperative architectures through case studies. The cases employed here are inspired by the Starlink constellation.
Finally, Sec.~\ref{sec:conclusion} concludes this paper with some remarks.

% \iffalse
% satellite constellation

% on orbit service is emerging

% northrop gramman example
% even more promising

% refuel is one application

% traditionally non coooperative 
% cite noncoop case
% coop case is relatively few
% cite noncoop

% to compare them, this paper

% the rest of the paper is organized as follows

% \fi

\section{Modeling Multi-Target On-Orbit Refueling Architectures}
\label{sec:problem_formulation}
% \subsection{Mathematical Model}
This section first introduces the considered on-orbit refuel mission architecture (Sec.~\ref{subsec:coopoverview}), followed by the derivation of the initial mass of a servicer (Sec.~\ref{subsec:imserv}) and total refuel mass (Sec.~\ref{subsec:mass_refuel}) in this architecture.
Then, Sec.~\ref{sec:analytical_comparison} analytically compares the derived mass between both cooperative and non-cooperative architectures.

\subsection{Overview of On-Orbit Refueling Architecture}
\label{subsec:coopoverview}
\begin{figure}[hbt!]
\center
\includegraphics[width=.9\textwidth]{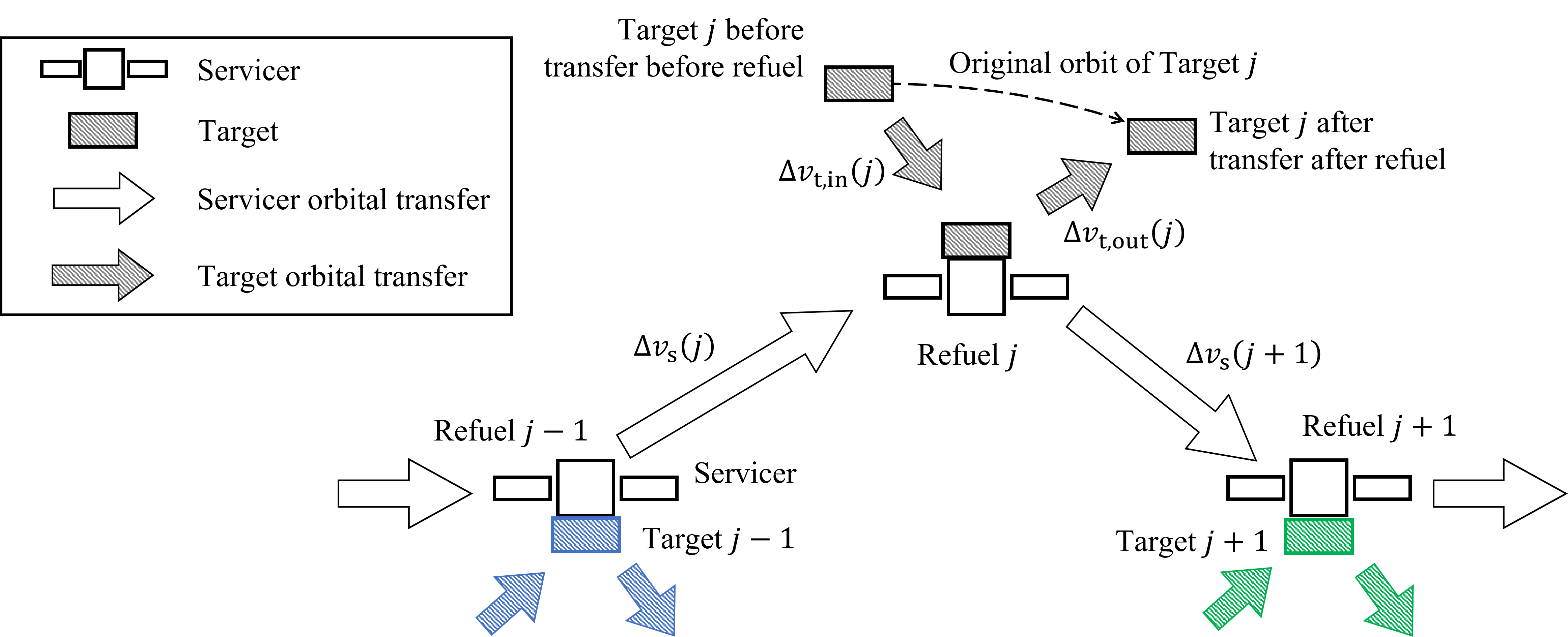}
\caption{Cooperative on-orbit refueling concept.}
\label{fig:coop_seq}
\end{figure}

% This paper assumes that target spacecraft are distributed in multiple orbital planes. 
A general multi-target on-orbit refueling architecture considered in this paper follows the  
% Based on \cite{zhao2017}, 
% This paper considers an on-orbit refueling following the 
% this subsection formulates the mathematical model of an on-orbit refuel mission that follows the 
sequence below:
\begin{enumerate}[label=\arabic*),ref=\arabic*]
  \item a single servicer is launched and transferred to the orbital plane of the 1st target;
  \item \label{seq:2}the servicer rendezvous with the target, refuel it, and move to the next target; 
  \item repeat Phase~\ref{seq:2}) until all targets are refueled;
  \item \label{seq:3}the servicer goes back to its original position with the original inclination and the relative argument of latitude.
\end{enumerate}
In the case of a cooperative refueling architecture, both the servicer and a target are capable of 
orbital transfers to rendezvous in Phase \ref{seq:2} (Fig. \ref{fig:coop_seq}), 
whereas in the conventional non-cooperative architecture, only the servicer can conduct orbital transfers to the rendezvous point.
Note that cooperative target spacecraft consume some fuel to fly to a rendezvous point and come back to their original position after refueling (See Fig. \ref{fig:coop_seq}), and therefore, they need to receive more fuel from the servicer than those in the non-cooperative architecture to compensate this extra fuel consumption.
% In general, 
% % for the cooperative architecture, 
% to receive the same amount of fuel at the end of the whole procedure, 
% each target in the cooperative architecture  needs to receive more fuel from the servicer than those in the non-cooperative architecture. This is because cooperative target spacecraft consume some fuel to fly to a rendezvous point and come back to their original position after refueling (See Fig. \ref{fig:coop_seq}).

% Note that the total refuel amount may differ between these two cases since targets need more fuel for their active orbital transfers for the cooperative case. 

\subsection{Initial Mass of Servicer}
\label{subsec:imserv}
To compare the cooperative and non-cooperative architectures introduced above, the mathematical expressions of the initial mass of the servicer in both architectures are derived in this subsection.
For the $j$\,th rendezvous (See Fig. \ref{fig:coop_seq}), the Tsiolkovsky rocket equation can be written as 
\begin{equation}
m_\mathrm{s,0}(j) =  m_\mathrm{s,f}(j) \exp\left(\frac{\Delta v_\mathrm{s} (j)}{I_\mathrm{sp, s} g_0}\right),
\label{eq:rocket}
\end{equation}
where $m_\mathrm{s,0}(j)$ and $m_\mathrm{s,f}(j)$ denote the mass of the servicer before and after the $j$\,th orbital transfer, respectively.
Since the initial servicer mass before the $j$\,th orbital transfer is the same as the servicer mass after the $\left(j-1\right)$\,th refuel, 
% they follow the following relationship,
\begin{equation}
m_\mathrm{s,0}(j) =  m_\mathrm{s,f}(j-1) - m_\mathrm{r}(j-1),
\label{eq:mass_conservation}
\end{equation}
where $m_\mathrm{r}(j-1)$ represents the total refuel mass that the servicer gave to the $\left(j-1\right)$\,th target.
Modifying Eq.~\eqref{eq:rocket} to the $\left(j-1\right)$\,th orbital transfer and combining it with Eq.~\eqref{eq:mass_conservation} leads to 
\begin{equation}
m_\mathrm{s,0}(j-1) =  \left(m_\mathrm{s,0}(j) + m_\mathrm{r}(j-1)\right) \exp\left(\frac{\Delta v_\mathrm{s} (j - 1)}{I_\mathrm{sp, s} g_0}\right).
\label{eq:mass_j-1_to_j}
\end{equation}
From this equation, assuming the number of targets is $n$, the initial mass of the servicer at the beginning of the whole sequence, $m_\mathrm{s,0}(1)$, is derived as follows:
\begin{equation}
m_\mathrm{s,0}(1) =  m_\mathrm{s,0}(n) \exp\left(\frac{\sum_{j = 1}^{n-1}\Delta v_\mathrm{s} (j)}{I_\mathrm{sp, s} g_0}\right) + \sum_{j = 1}^{n-1} \left( m_\mathrm{r}(j) \exp\left(\frac{\sum_{k = 1}^{j}\Delta v_\mathrm{s} (k)}{I_\mathrm{sp, s} g_0}\right) \right).
\label{eq:mass_i_to_1}
\end{equation}
We further assume the final mass of the servicer is given (e.g., the dry mass) and treat Phase~\ref{seq:3}) as the $(n+1)$\,th orbital transfer. Thus, the given final mass is \,$m_\mathrm{s,F} = m_\mathrm{s,f}(n+1)$. From Eqs.~\eqref{eq:rocket} and \eqref{eq:mass_i_to_1}, the relationship between the initial mass,\,$m_\mathrm{s,I} = m_\mathrm{s,0}(1)$, and the final mass,  of the servicer of the whole scenario\,$m_\mathrm{s,F}$, is derived as follows:
\begin{align}
% m_\mathrm{s,0}(1) &=  m_\mathrm{s,0}(n+1) \exp\left(\frac{\sum_{j = 1}^{n}\Delta v_\mathrm{s} (j)}{I_\mathrm{sp, s} g_0}\right) + \sum_{j = 1}^{n} \left( m_\mathrm{r}(j) \exp\left(\frac{\sum_{k = 1}^{j}\Delta v_\mathrm{s} (k)}{I_\mathrm{sp, s} g_0}\right) \right),\nonumber\\
% \therefore 
m_\mathrm{s,I} = m_\mathrm{s,F} \exp\left(\frac{\sum_{j = 1}^{n+1}\Delta v_\mathrm{s} (j)}{I_\mathrm{sp, s} g_0}\right) + \sum_{j = 1}^{n} \left( m_\mathrm{r}(j) \exp\left(\frac{\sum_{k = 1}^{j}\Delta v_\mathrm{s} (k)}{I_\mathrm{sp, s} g_0}\right) \right).
\label{eq:mass_final_to_1}
\end{align}

\subsection{Total Refuel Mass}
\label{subsec:mass_refuel}
\begin{figure}[bt!]
\center
\includegraphics[width=.7\textwidth]{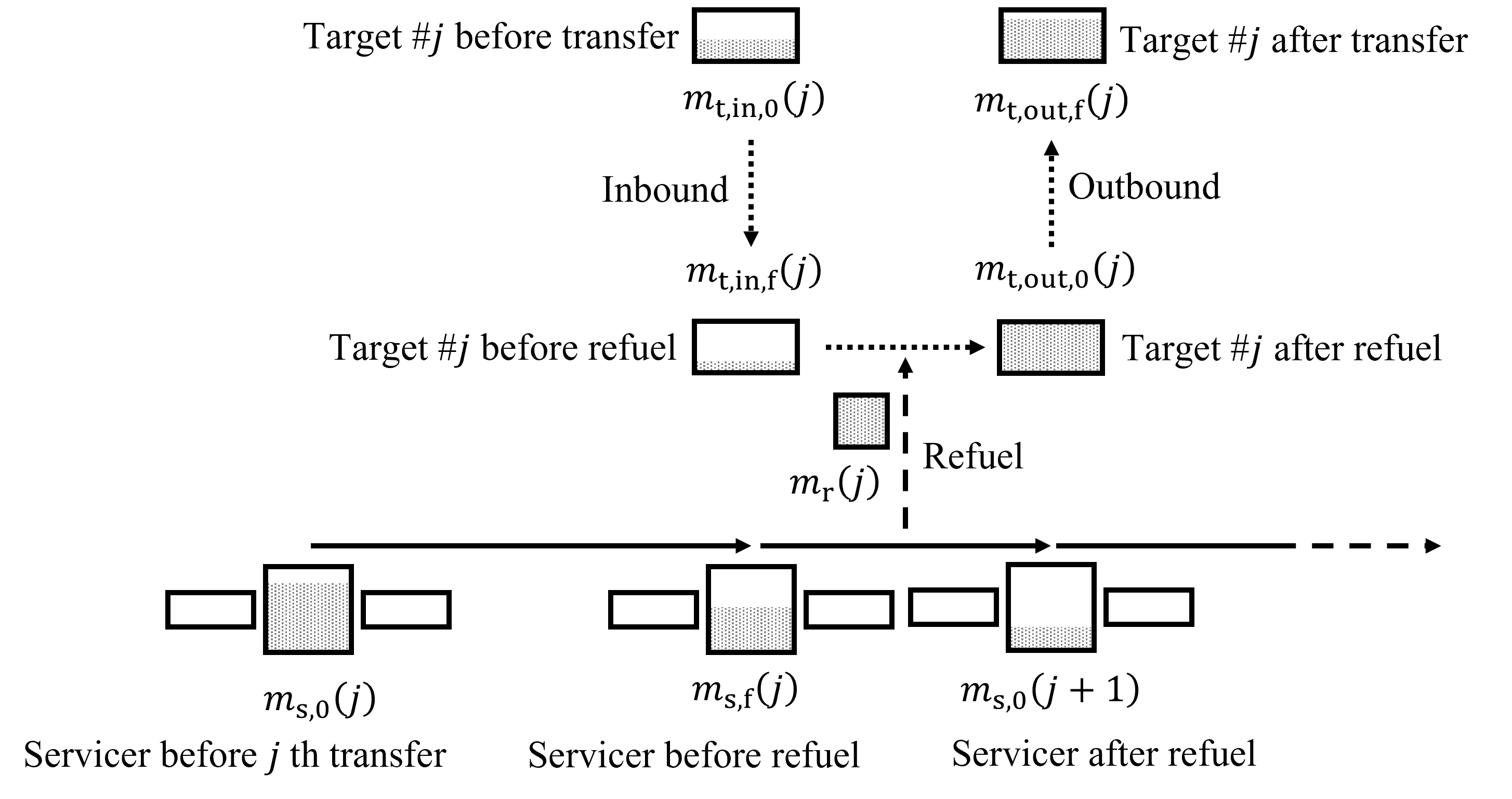}
\caption{Change in mass during the $j$\,th cooperative refuel.}
\label{fig:coop_seq_mass}
\end{figure}

For the cooperative architecture, a target comes from its original position to a rendezvous spot, and then the target goes back to the position where it would be if there was no active orbital transfer (see Fig. \ref{fig:coop_seq}). The rocket equations of both transfers by Target $j$ can be written as
\begin{align}
\label{eq:rocket_target}
m_\mathrm{t,in,f}(j) &=  m_\mathrm{t,in,0}(j)  \exp\left(\frac{-\Delta v_\mathrm{t, in} (j)}{I_\mathrm{sp, t} g_0}\right),\\
m_\mathrm{t,out,f}(j) &=  m_\mathrm{t,out,0}(j) \exp\left(\frac{-\Delta v_\mathrm{t, out} (j)}{I_\mathrm{sp, t} g_0}\right),
\end{align}
where $m_\mathrm{t,in,f}(j)$ and $m_\mathrm{t,in,0}(j)$ are the initial and the final mass of the first (inbound) transfer, $m_\mathrm{t,out,f}(j)$ and $m_\mathrm{t,in,0}(j)$ are those of the second (outbound) transfer (Fig. \ref{fig:coop_seq_mass}), and $\Delta v_\mathrm{t, in}$ and $\Delta v_\mathrm{t, out}$ are the $\Delta v$ required for each transfer. To fairly compare the cooperative and the non-cooperative architectures, the final mass of Target $j$ (after the second transfer),\,$m_\mathrm{t,out,f}(j)$, is assumed to be the same as the post-refuel mass of that Target in the non-cooperative architecture, i.e., a target in either architecture have the same amount of ``usable'' fuel supplied after every orbital transfer is completed.
Therefore, the relationship between the initial mass of the Target $j$, $m_\mathrm{t,in,0}(j)$, and the final mass of the Target $j$, $m_\mathrm{t,out,f}(j)$, becomes
\begin{equation}
m_\mathrm{t,out,f}(j) = m_\mathrm{t,in,0}(j) + m_\mathrm{req}(j),
\label{eq:mass_target_final}
\end{equation}
where $m_\mathrm{req}(j)$ is the required ``usable'' refuel amount by Target $j$ (i.e., the total refuel amount for the non-cooperative architecture). 
Furthermore, the relationship between the final mass after the inbound transfer and the initial mass before the outbound transfer is expressed as 
\begin{equation}
m_\mathrm{t,out,0}(j) = m_\mathrm{t,in,f}(j) + m_\mathrm{r,c}(j).
\label{eq:mass_t_refuel}
\end{equation}
where $m_\mathrm{r,c}(j)$ is $m_\mathrm{r}(j)$ for the cooperative architecture. From Eqs.~\eqref{eq:rocket_target} - \eqref{eq:mass_t_refuel}, the total refuel mass for Target $j$ is derived as follows:
\begin{equation}
\begin{split}
% m_\mathrm{t,in,0}(j) + m_\mathrm{req}(j) &= \left(m_\mathrm{t,in,0}(j)  \exp\left(\frac{-\Delta v_\mathrm{t, in} (j)}{I_\mathrm{sp, t} g_0}\right) + m_\mathrm{r}(j)\right)\exp\left(\frac{-\Delta v_\mathrm{t, out} (j)}{I_\mathrm{sp, t} g_0}\right),\\
% \therefore 
m_\mathrm{r}(j) &= \left( m_\mathrm{t,in,0}(j) + m_\mathrm{req}(j) \right) \exp\left(\frac{\Delta v_\mathrm{t, out} (j)}{I_\mathrm{sp, t} g_0}\right) - m_\mathrm{t,in,0}(j)  \exp\left(\frac{-\Delta v_\mathrm{t, in} (j)}{I_\mathrm{sp, t} g_0}\right).
\end{split}
\label{eq:mass_total_refuel}
\end{equation}
Substituting Eq.~\eqref{eq:mass_total_refuel} into Eq.~\eqref{eq:mass_final_to_1},
the initial mass of the servicer in the cooperative architecture, $m_\mathrm{s,I,c}$, can be analytically derived:
\begin{align}
\nonumber
m_\mathrm{s,I,c}& = m_\mathrm{s,F} \exp\left(\frac{\sum_{j = 1}^{n+1}\Delta v_\mathrm{s} (j)}{I_\mathrm{sp, s} g_0}\right)\\ 
&+ \sum_{j = 1}^{n} \left( \left[ \left( m_\mathrm{t,in,0}(j) + m_\mathrm{req}(j) \right) \exp\left(\frac{\Delta v_\mathrm{t, out} (j)}{I_\mathrm{sp, t} g_0}\right) - m_\mathrm{t,in,0}(j)  \exp\left(\frac{-\Delta v_\mathrm{t, in} (j)}{I_\mathrm{sp, t} g_0}\right)\right] \exp\left(\frac{\sum_{k = 1}^{j}\Delta v_\mathrm{s} (k)}{I_\mathrm{sp, s} g_0}\right) \right).
\label{eq:mass_final_coop}
\end{align}

On the other hand, for the non-cooperative architecture, 
\begin{equation}
m_\mathrm{r,n}(j) =  m_\mathrm{req}(j).
\label{eq:mass_total_refuel_noncoop}
\end{equation}
where $m_\mathrm{r,n}(j)$ is $m_\mathrm{r}(j)$ for the non-cooperative architecture.
Therefore, substituting Eq.~\eqref{eq:mass_total_refuel_noncoop} into Eq.~\eqref{eq:mass_final_to_1},
the initial mass of the servicer in the conventional non-cooperative architecture, $m_\mathrm{s,I,n}$, can be analytically derived:
\begin{equation}
    m_\mathrm{s,I,n}= m_\mathrm{s,F} \exp\left(\frac{\sum_{j = 1}^{n+1}\Delta v_\mathrm{s,n} (j)}{I_\mathrm{sp, s} g_0}\right) + \sum_{j = 1}^{n} \left( m_\mathrm{req}(j) \exp\left(\frac{\sum_{k = 1}^{j}\Delta v_\mathrm{s,n} (k)}{I_\mathrm{sp, s} g_0}\right) \right).
\label{eq:mass_final_noncoop}
\end{equation}

% \begin{align}
% \nonumber
% m_\mathrm{s,I,c}& = m_\mathrm{s,F} \exp\left(\frac{\sum_{j = 1}^{n+1}\Delta v_\mathrm{s} (j)}{I_\mathrm{sp, s} g_0}\right)\\ 
% &+ \sum_{j = 1}^{n} \left( \left[ \left( m_\mathrm{t,in,0}(j) + m_\mathrm{req}(j) \right) \exp\left(\frac{\Delta v_\mathrm{t, out} (j)}{I_\mathrm{sp, t} g_0}\right) - m_\mathrm{t,in,0}(j)  \exp\left(\frac{-\Delta v_\mathrm{t, in} (j)}{I_\mathrm{sp, t} g_0}\right)\right] \exp\left(\frac{\sum_{k = 1}^{j}\Delta v_\mathrm{s} (k)}{I_\mathrm{sp, s} g_0}\right) \right).
% \label{eq:mass_final_coop}
% \end{align}
% For the non-cooperative architecture, from Eqs.~\eqref{eq:mass_final_to_1} and \eqref{eq:mass_total_refuel_noncoop}, 
% the initial mass of the servicer,\,$m_\mathrm{s,I,n}$, is expressed as 
% \begin{equation}
%     m_\mathrm{s,I,n}= m_\mathrm{s,F} \exp\left(\frac{\sum_{j = 1}^{n+1}\Delta v_\mathrm{s,n} (j)}{I_\mathrm{sp, s} g_0}\right) + \sum_{j = 1}^{n} \left( m_\mathrm{req}(j) \exp\left(\frac{\sum_{k = 1}^{j}\Delta v_\mathrm{s,n} (k)}{I_\mathrm{sp, s} g_0}\right) \right).
% \label{eq:mass_final_noncoop}
% \end{equation}

\subsection{Comparison of Initial Mass of Servicer}
\label{sec:analytical_comparison}
% This section analytically compares the initial mass of the servicer between both cooperative and non-cooperative architectures.
% From Eqs.~\eqref{eq:mass_final_to_1} and \eqref{eq:mass_total_refuel_noncoop}, the initial mass of the servicer in the non-cooperative architecture,\,$m_\mathrm{s,I,n}$, is expressed as 
% \begin{equation}
%     m_\mathrm{s,I,n}= m_\mathrm{s,F} \exp\left(\frac{\sum_{j = 1}^{n+1}\Delta v_\mathrm{s,n} (j)}{I_\mathrm{sp, s} g_0}\right) + \sum_{j = 1}^{n} \left( m_\mathrm{req}(j) \exp\left(\frac{\sum_{k = 1}^{j}\Delta v_\mathrm{s,n} (k)}{I_\mathrm{sp, s} g_0}\right) \right).
% \label{eq:mass_final_noncoop}
% \end{equation}
For designing a refueling architecture, a key cost metric would be the initial mass (wet mass) of the servicer. To obtain conditions that the cooperative case is favored, the following inequality is considered:
\begin{equation}
    m_\mathrm{s,I,c}\le m_\mathrm{s,I,n}.
\end{equation}
% \begin{multline}
%     m_\mathrm{s,I,c}<m_\mathrm{s,I,n}\Leftrightarrow m_\mathrm{s,F} \exp\left(\frac{\sum_{j = 1}^{n+1}\Delta v_\mathrm{s,c} (j)}{I_\mathrm{sp, s} g_0}\right) + \sum_{j = 1}^{n} \left( m_\mathrm{r,c}(j) \exp\left(\frac{\sum_{k = 1}^{j}\Delta v_\mathrm{s,c} (k)}{I_\mathrm{sp, s} g_0}\right) \right)\\
%     < m_\mathrm{s,F} \exp\left(\frac{\sum_{j = 1}^{n+1}\Delta v_\mathrm{s,n} (j)}{I_\mathrm{sp, s} g_0}\right) + \sum_{j = 1}^{n} \left( m_\mathrm{req}(j) \exp\left(\frac{\sum_{k = 1}^{j}\Delta v_\mathrm{s,n} (k)}{I_\mathrm{sp, s} g_0}\right) \right).
% \label{eq:mass_ineq1}
% \end{multline}
From Eqs.~\eqref{eq:mass_final_coop} and \eqref{eq:mass_final_noncoop}, 
% From Eqs.~\eqref{eq:mass_final_to_1}, \eqref{eq:mass_total_refuel},  and \eqref{eq:mass_total_refuel_noncoop},
assuming all targets have the same initial mass (i.e., $m_\mathrm{t,in,0}(1) = m_\mathrm{t,in,0}(2) = \cdots = m_\mathrm{t,I}$), and they require the same amount of refuel (i.e., $m_\mathrm{req}(1) = m_\mathrm{req}(2) = \cdots = m_\mathrm{req}$), the modified version of this inequality is expressed as follows:
% {\color{red} add explanation HERE!!!!!!!.}
\begin{align}
    % \begin{split}
    %     &{}m_\mathrm{s,F} \exp\left(\frac{\sum_{j = 1}^{n+1}\Delta v_\mathrm{s,c} (j)}{I_\mathrm{sp, s} g_0}\right) + \sum_{j = 1}^{n} \left\{ \left[\left( m_\mathrm{t,I} + m_\mathrm{req} \right) \exp\left(\frac{\Delta v_\mathrm{t, out} (j)}{I_\mathrm{sp, t} g_0}\right) - m_\mathrm{t,I}  \exp\left(\frac{-\Delta v_\mathrm{t, in} (j)}{I_\mathrm{sp, t} g_0}\right)\right] \exp\left(\frac{\sum_{k = 1}^{j}\Delta v_\mathrm{s,c} (k)}{I_\mathrm{sp, s} g_0}\right) \right\}\nonumber\\
    %     &{}\le   m_\mathrm{s,F} \exp\left(\frac{\sum_{j = 1}^{n+1}\Delta v_\mathrm{s,n} (j)}{I_\mathrm{sp, s} g_0}\right) + \sum_{j = 1}^{n} \left(m_\mathrm{req} \exp\left(\frac{\sum_{k = 1}^{j}\Delta v_\mathrm{s,n} (k)}{I_\mathrm{sp, s} g_0}\right) \right),\nonumber
    % \end{split}\\
    \begin{split}
    \label{eq:mass_ineq2}
        &{}\frac{m_\mathrm{s,F}}{m_\mathrm{t,I}} \left(\exp\left(\frac{\sum_{j = 1}^{n+1}\Delta v_\mathrm{s,c} (j)}{I_\mathrm{sp, s} g_0}\right) - \exp\left(\frac{\sum_{j = 1}^{n+1}\Delta v_\mathrm{s,n} (j)}{I_\mathrm{sp, s} g_0}\right)\right)\\
        &{}\le \sum_{j = 1}^{n} \left\{ \left( \frac{m_\mathrm{req}}{m_\mathrm{t,I}} \exp\left(\frac{\sum_{k = 1}^{j}\Delta v_\mathrm{s,n} (k)}{I_\mathrm{sp, s} g_0}\right) \right) -   \left[\left(1 +\frac{m_\mathrm{req}}{m_\mathrm{t,I}}  \right) \exp\left(\frac{\Delta v_\mathrm{t, out} (j)}{I_\mathrm{sp, t} g_0}\right) - \exp\left(\frac{-\Delta v_\mathrm{t, in} (j)}{I_\mathrm{sp, t} g_0}\right)\right] \exp\left(\frac{\sum_{k = 1}^{j}\Delta v_\mathrm{s,c} (k)}{I_\mathrm{sp, s} g_0}\right)\right\} .
    \end{split}
\end{align}

If Inequality \eqref{eq:mass_ineq2} is satisfied, the initial mass of the servicer can be smaller for the cooperative case than that of the conventional non-cooperative mission architecture. Note that this inequality depends on the required $\Delta v$ of both the servicer, and the targets, as well as the final mass of the servicer, the initial mass of the targets, the required refuel mass, and the specific impulse of both the servicer and the targets.

Inequality \eqref{eq:mass_ineq2} also indicates that one critical parameter to evaluate the architectures is the ratio between the final servicer mass and the initial target mass,~$m_\mathrm{s,F}/m_\mathrm{t,I}$. The numerator of this ratio, the final servicer mass, is the dry mass of the servicer, whereas the denominator of this ratio, the initial target mass, is the target mass before the entire fuel campaign starts. Note that, despite the name ``initial,'' the initial target mass is not a fully-fueled target; rather it is the mass of the target when refueling is requested, at which point its fuel tank can be near empty although not necessarily completely empty, especially for the cooperative architecture. 

From Inequality \eqref{eq:mass_ineq2}, we can find a crossover point of this mass ratio,~$m_\mathrm{s,F}/m_\mathrm{t,I}$ beyond which the cooperative architecture leads to a smaller initial servicer mass than the non-cooperative architecture. In this paper, we call this value the ``critical mass ratio,'' denoted by $\alpha$. 
Refer to Appendix\ref{app:crit_mass_ratio} for the analytical expression of the critical mass ratio.
In general, the required $\Delta v$ for the servicer is lower in the cooperative architectures. Therefore, from Inequality \eqref{eq:mass_ineq2}, if $m_\mathrm{s,F}/m_\mathrm{t,I}$ is smaller than the critical mass ratio (i.e., a light-weight servicer and heavy-weight targets), the initial servicer mass in the non-cooperative architecture can be smaller than the cooperative architecture, i.e., the non-cooperative architecture is favored.
Instead, if the value is bigger than the critical mass ratio (i.e., a heavy-weight servicer and light-weight targets), the cooperative architecture can be beneficial.
Therefore, calculating the ratio $m_\mathrm{s,F}/m_\mathrm{t,I}$ and comparing it with the critical mass ratio $\alpha$ at the phase of conceptual study or preliminary analysis of designing an on-orbit refueling mission gives mission architects a recommended architecture between these two mission architectures.
This ratio, $m_\mathrm{s,F}/m_\mathrm{t,I}$, can be estimated from the mass of targets requested from customers and the required payloads of the servicer such as a robot arm for capturing the targets.

% if the ratio between the final servicer mass and the initial target mass,~$m_\mathrm{s,F}/m_\mathrm{t,I}$, satisfies 
% {\color{red} add explanation of critical mass ratio HERE!!!!!!!.}
% \begin{equation}
% \label{eq:mass_critical}
%     \frac{m_\mathrm{s,F}}{m_\mathrm{t,I}} =  \frac{ \sum_{j = 1}^{n} \left\{ \left( \frac{m_\mathrm{req}}{m_\mathrm{t,I}} \exp\left(\frac{\sum_{k = 1}^{j}\Delta v_\mathrm{s,n} (k)}{I_\mathrm{sp, s} g_0}\right) \right) -   \left[\left(1 +\frac{m_\mathrm{req}}{m_\mathrm{t,I}}  \right) \exp\left(\frac{\Delta v_\mathrm{t, out} (j)}{I_\mathrm{sp, t} g_0}\right) - \exp\left(\frac{-\Delta v_\mathrm{t, in} (j)}{I_\mathrm{sp, t} g_0}\right)\right] \exp\left(\frac{\sum_{k = 1}^{j}\Delta v_\mathrm{s,c} (k)}{I_\mathrm{sp, s} g_0}\right)\right\}  }{\exp\left(\frac{\sum_{j = 1}^{n+1}\Delta v_\mathrm{s,c} (j)}{I_\mathrm{sp, s} g_0}\right) - \exp\left(\frac{\sum_{j = 1}^{n+1}\Delta v_\mathrm{s,n} (j)}{I_\mathrm{sp, s} g_0}\right)},
% \end{equation}
% both architectures require the same initial servicer mass. In this paper, we call this value the "critical mass ratio." 
% In general, the required $\Delta v$ for the servicer is lower in the cooperative architectures. Therefore, from Inequality \eqref{eq:mass_ineq2}, if $m_\mathrm{s,F}/m_\mathrm{t,I}$ is smaller than the critical mass ratio, the initial servicer mass in the non-cooperative architecture can be smaller than the cooperative architecture, i.e., the non-cooperative architecture is favored.
% Instead, if the value is bigger than the critical mass ratio, the cooperative architecture can be benefitial.

\section{Case Study}
\label{sec:case}
In this section, case studies of on-orbit refueling missions in multiple architectures are investigated.
First, Sec.~\ref{subsec:case_overview} introduces the overview of a target LEO satellite constellation and mission architectures considered in this paper.
The initial mass of the servicer in each architecture is summarized in Sec.~\ref{subsec:case_init_mass}.
Finally, the sensitivity of each parameter to the the critical mass ratio is discussed in Sec.~\ref{subsec:case_sensitivity}.

\subsection{Case Study Overview}
\label{subsec:case_overview}
\begin{figure}[hbt!]
\center
\includegraphics[width=.9\textwidth]{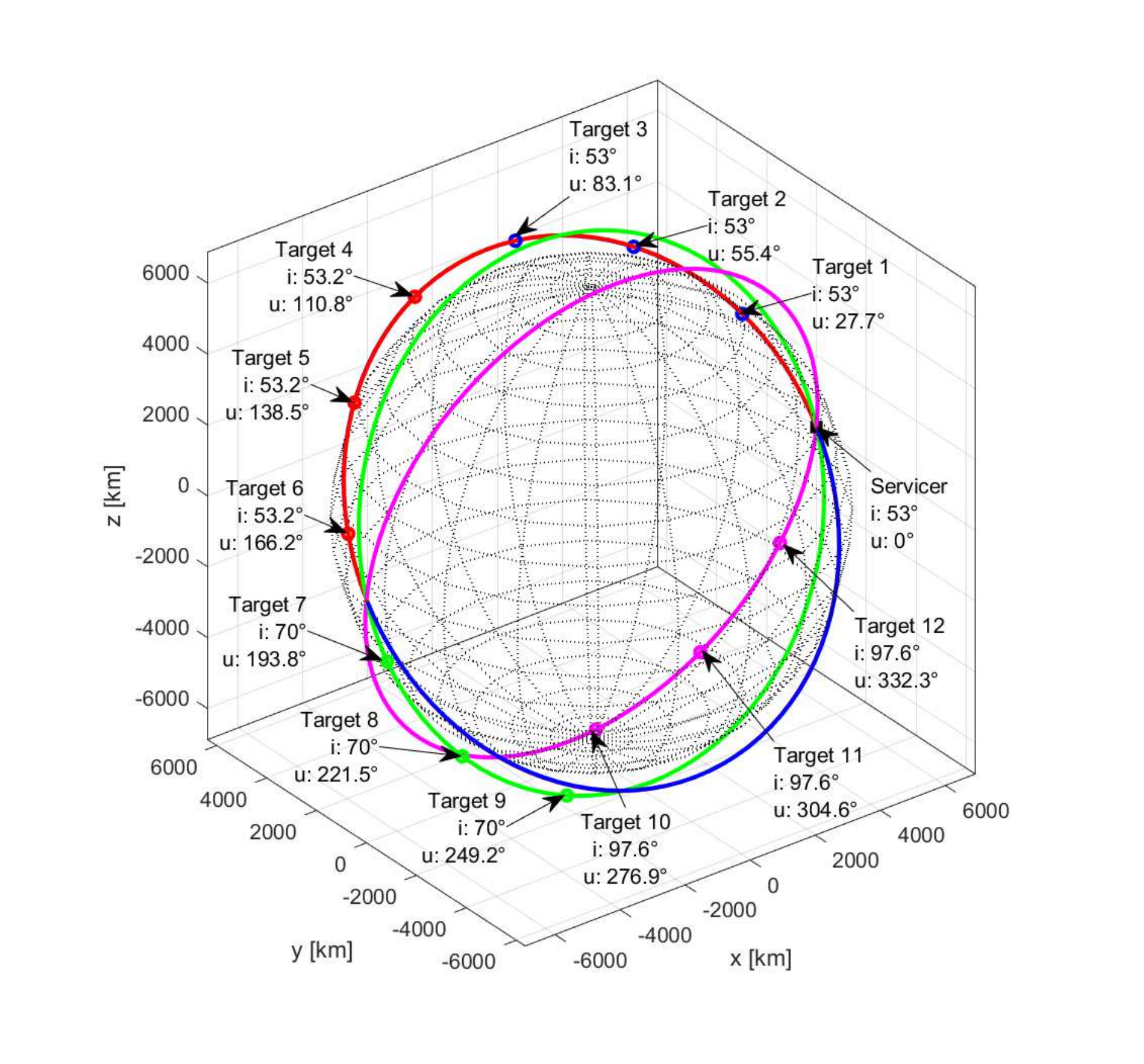}
\caption{Initial positions of the targets and the servicer in the Earth-Centered Inertial (ECI) frame.}
\label{fig:target_position}
\end{figure}

To compare cooperative and conventional non-cooperative on-orbit refueling architectures,  simple multi-target refueling missions in LEOs are considered in this paper.
% Each mission has a different number of targets (e.g., $n = $6, 9, or 12) in a LEO constellation.
% with a single servicer and a different number of targets (i.e., $n = $6, 9, and 12) are considered in this paper.
The target spacecraft are assumed to be deployed in multiple orbital planes with different inclination angles (i.e., $i=$53°, 53.2°, 70°, or 97.6° inspired by Starlink satellites).
% To compare the architectures simply,
We assume that
% For the sake of simplicity, 
all of these orbits are circular with an altitude of 550 km, and the targets share the same Keplerian orbital elements except for the inclination and the argument of latitude. 
Table \ref{tab:parameters} lists the parameters used in this case study, and Fig.~\ref{fig:target_position} illustrates the initial positions of the targets and the servicer considered in the Earth-Centered Inertial (ECI) frame.
In this case study, the servicer separately refuels to Targets 1 to $n$ in Fig.~\ref{fig:target_position}.

\begin{table}[hbt!]
\center
\caption{Parameters employed for the case study}
\label{tab:parameters}
\begin{tabular}{lcr}
\thickhline
\multicolumn{2}{c}{Parameter}&Assumed value\\\hline
\multicolumn{3}{c}{\textit{Constellation Parameters}}\\
number of targets        &$n$     & 1 - 12\\
initial inclination [deg] &$i$& 53, 53.2, 70, 97.6\\
initial argument of latitude  &$u$&  See Fig.~\ref{fig:target_position} \\
altitude [km] && 550\\
\multicolumn{3}{c}{\textit{Refueling Mission Parameters}}\\
specific impulse of target [s] & $I_\mathrm{sp, t}$ & 300\\
specific impulse of servicer [s] & $I_\mathrm{sp, s}$ & 300\\
initial mass of each target [kg] & $m_\mathrm{t, I}$ & 1000\\
required refuel amount [kg] & $m_\mathrm{req}$ & 200\\\hline
\end{tabular}
\end{table}

To accomplish refueling to one target, either or both of the target and the servicer need to move to a rendezvous position. In this paper, we assume that these orbital transfers consist of two phases: inclination change and coplanar phasing.
For circular orbits, the required $\Delta v$ to change the inclination,\,$\Delta v_\mathrm{inc}$, is calculated as
\begin{equation}
\Delta v_\mathrm{inc} = 2 v \sin{\left( \frac{\Delta i}{2}\right)},
\label{eq:deltav_i}
\end{equation}
where $v$ and $\Delta i$ denote the velocity of the circular orbit $\sqrt{\mu/r}$ and the change in inclination, respectively.
The $\Delta v$ required for phasing maneuvers,\,$\Delta v_\mathrm{pha}$, is 
\begin{equation}
\Delta v_\mathrm{pha} = 2 \left | \sqrt{\frac{\mu}{r}} - \sqrt{\mu \left(\frac{2}{r} - \frac{1}{a}\right )} \right  |,
\label{eq:deltav_pha}
\end{equation}
where
\begin{align}
a &= \left( \frac{\alpha + 2 \pi k_2}{2 \pi k_1}\right) ^ {2/3} r,\\
\alpha &= 2\pi - \Delta u
\label{eq:deltav_pha_a}
\end{align}
$k_1$ and $k_2$ denote the number of complete revolutions during the phasing by the servicer and the target, respectively, and $\Delta u \in [0, 2\pi)$ denote the initial difference in the argument of latitude measured from the one that chases the other one. For the non-cooperative architecture, the servicer is always the one that chases the target, and in a cooperative architecture, it depends on their transfer strategies.
Please refer to Refs.~\cite{vall2013,sart2021} for the detailed derivation of Eq.~\eqref{eq:deltav_pha}.

\begin{figure}[hbt!]
\center
\includegraphics[width=.99\textwidth]{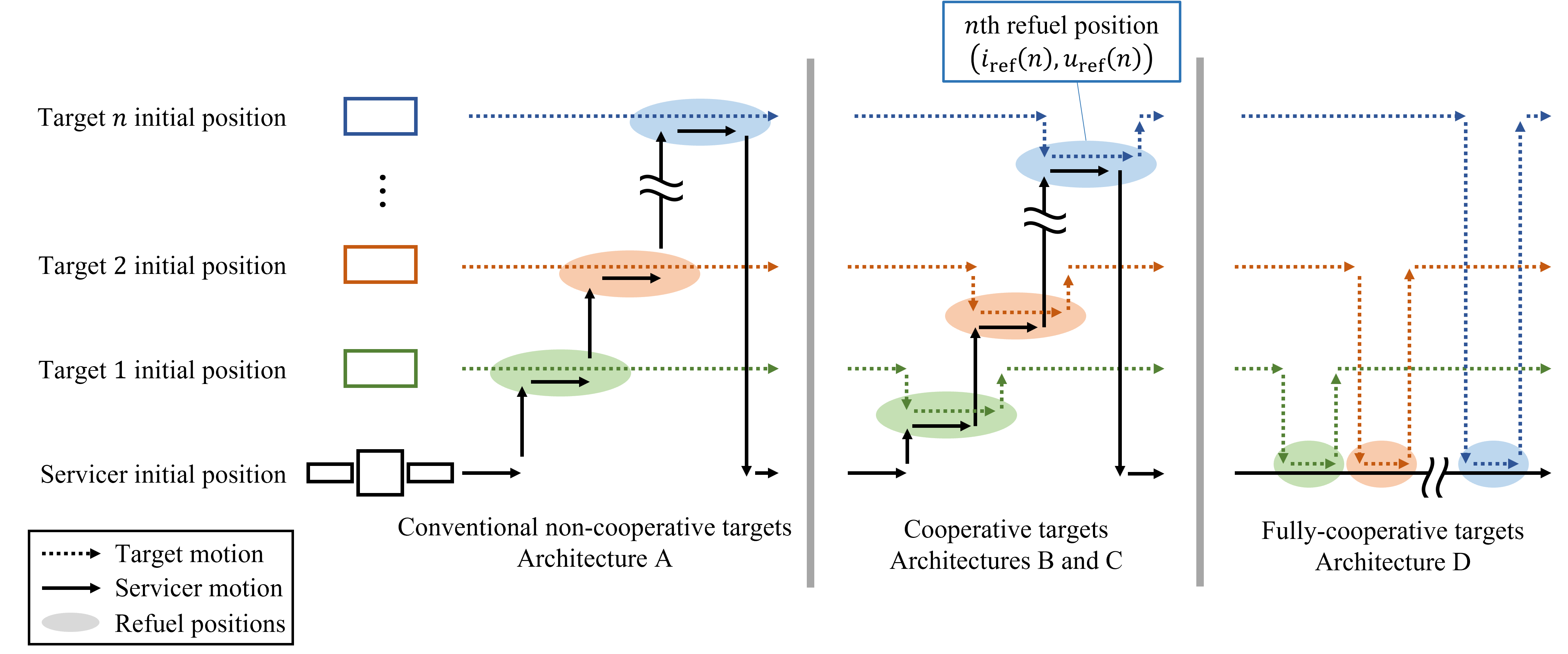}
\caption{Mission architectures of the multi-target refueling considered in this paper.}
\label{fig:architecture}
\end{figure}

This paper considers the following refueling architectures for the given LEO constellation:
\begin{enumerate}[label=\Alph*),ref=\Alph*]
\item \label{arc:trad}the conventional non-cooperative architecture where only the servicer actively conducts the orbital transfers;
\item a cooperative architecture where the inclination change is done by the targets and the phasing is done by the servicer;
\item a cooperative architecture where the inclination change is done by the servicer and the phasing is done by the targets;
\item \label{arc:extreme}an architecture where all targets come to the servicer while the servicer remains in the same orbit (i.e., fully cooperative targets and passive servicer ``depot''); and
\item \label{arc:optimum}the optimal cooperative case that minimizes the initial mass of the servicer.
\end{enumerate}
Figure \ref{fig:architecture} visually summarizes Architectures \ref{arc:trad}-\ref{arc:extreme}.
Architecture \ref{arc:optimum} can be any of Architectures \ref{arc:trad}-\ref{arc:extreme} or an architecture where both the servicer and the targets conduct both the inclination change and phasing maneuvers if that can reduce the initial mass of the servicer further. Note that the order of the targets to be refueled is the same among all architectures. The consideration of the optimum of this order is left for future work.

As mentioned above, we consider both the inclination change and the coplanar phasing for each orbital transfer for refueling. The inclination $i$ and the relative argument of latitude $u$ to the initial argument of latitude of the servicer where the $j$\,th refuel is performed are denoted by $i_\mathrm{ref}(j)$ and $u_\mathrm{ref}(j)$, respectively, $i_\mathrm{init,t}(j)$ and $u_\mathrm{init,t}(j)$ represent the initial inclination and the initial relative argument of latitude of Target $j$, and $i_\mathrm{init,s}$ and $u_\mathrm{init,s}$ are those of the servicer. With this notation, each architecture has the following conditions:
\begin{enumerate}[label=\Alph*),ref=\Alph*]
\item $\forall j, i_\mathrm{ref}(j) = i_\mathrm{init,t}(j), u_\mathrm{ref}(j) = u_\mathrm{init,t}(j)$,
\item $\forall j, i_\mathrm{ref}(j) = i_\mathrm{init,s}, u_\mathrm{ref}(j) = u_\mathrm{init,t}(j)$,
\item $\forall j, i_\mathrm{ref}(j) = i_\mathrm{init,t}(j), u_\mathrm{ref}(j) = u_\mathrm{init,s}$, and
\item $\forall j, i_\mathrm{ref}(j) = i_\mathrm{init,s}, u_\mathrm{ref}(j) = u_\mathrm{init,s}$.
\end{enumerate}
For Architecture \ref{arc:optimum}, $i_\mathrm{ref}(j)$ and $u_\mathrm{ref}(j)$
are the design variables for the optimization, i.e., there are $2n$ design variables to be optimized. To find the optimal rendezvous points expressed by these parameters that minimize the initial mass of the servicer,\,$m_\mathrm{s,I}$, we use MATLAB's \texttt{MultiStart} solver.

\subsection{Results}
\label{subsec:case_init_mass}

\subsubsection{Initial Mass of Servicer}

We first examine the initial mass of the servicer and particularly how it is impacted by the ratio between the final servicer mass (i.e., servicer dry mass) and the initial target mass (i.e., target mass before the entire refuel campaign),~$m_\mathrm{s,F}/m_\mathrm{t,I}$.
The initial mass of the servicer is a summation of the following: dry mass, required refuel amount defined in this mission (i.e., $nm_\mathrm{req}$), consumed fuel mass by the servicer's orbital transfer $m_\mathrm{con, s}$ and the summation of consumed fuel mass by the targets' orbital transfer $m_\mathrm{con, t}$. Note that, for Architecture A, $m_\mathrm{con, t} = 0$, and for Architecture D, $m_\mathrm{con, s} = 0$.
Since the drymass and $nm_\mathrm{req}$ are fixed values between all architectures introduced in Sec.~\ref{subsec:case_overview}, we focused on the summation of $m_\mathrm{con, s}$ and $m_\mathrm{con, t}$, the ``variable fuel mass.''

\begin{figure}[hbt!]
\begin{subfigure}{.49\textwidth}
\centering
\includegraphics[width=.99\textwidth]{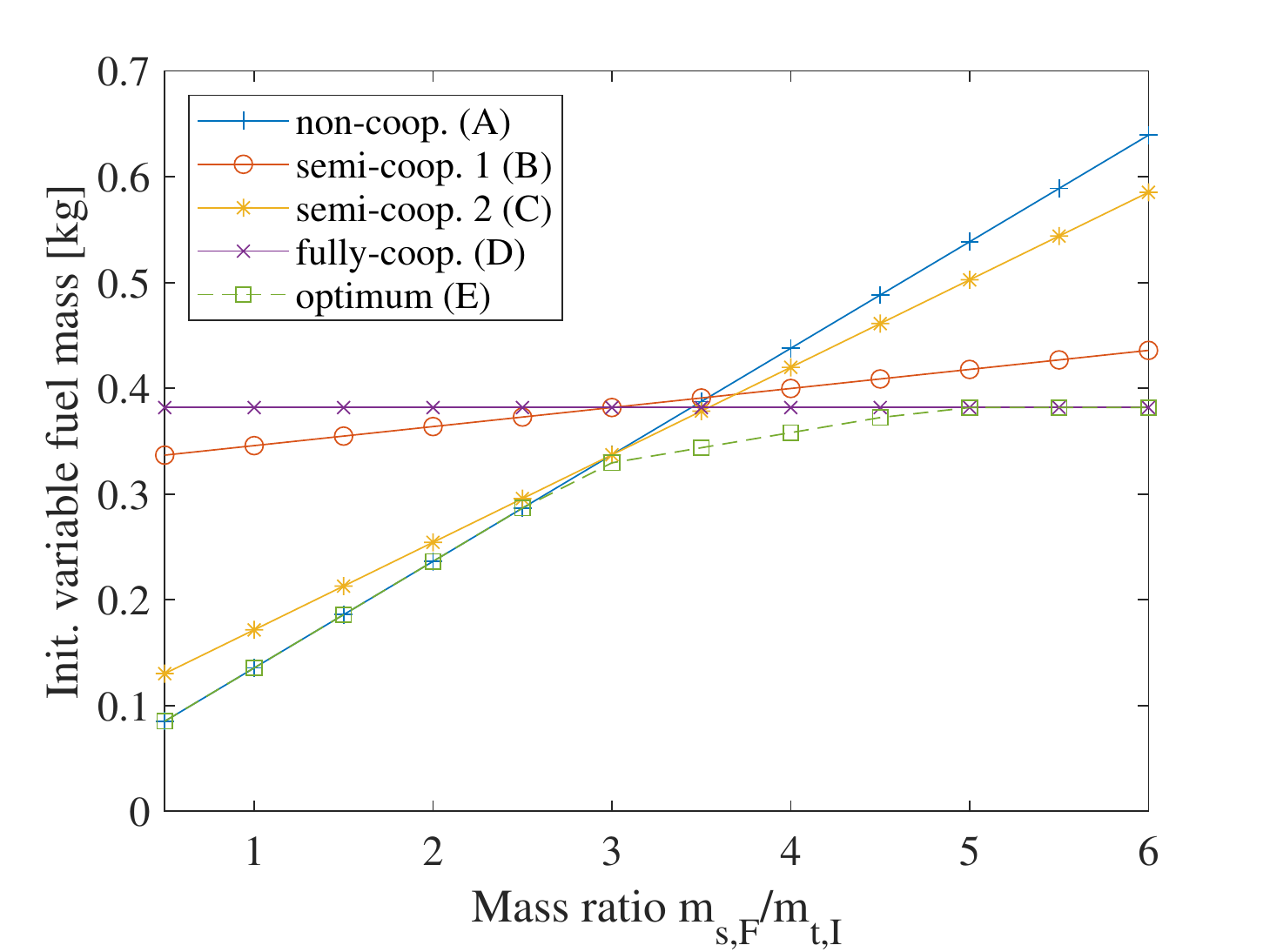}
\caption{$n=6$.}
\label{fig:case_fuelmass_abs_6}
\end{subfigure}
\begin{subfigure}{.49\textwidth}
\centering
\includegraphics[width=.99\textwidth]{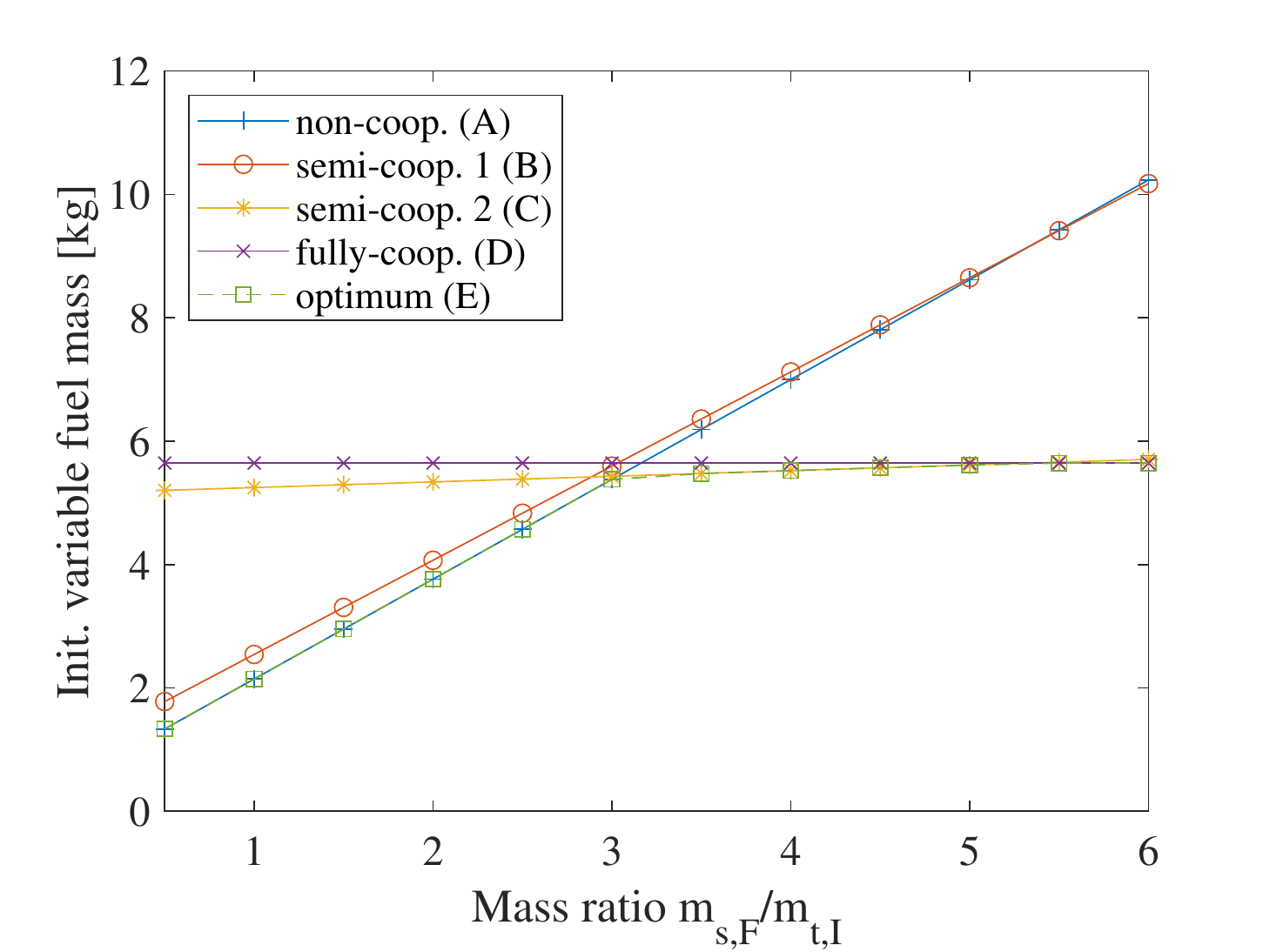}
\caption{$n=9$.}
\label{fig:case_fuelmass_abs_9}
\end{subfigure}\\
\begin{subfigure}{.49\textwidth}
\centering
\includegraphics[width=.99\textwidth]{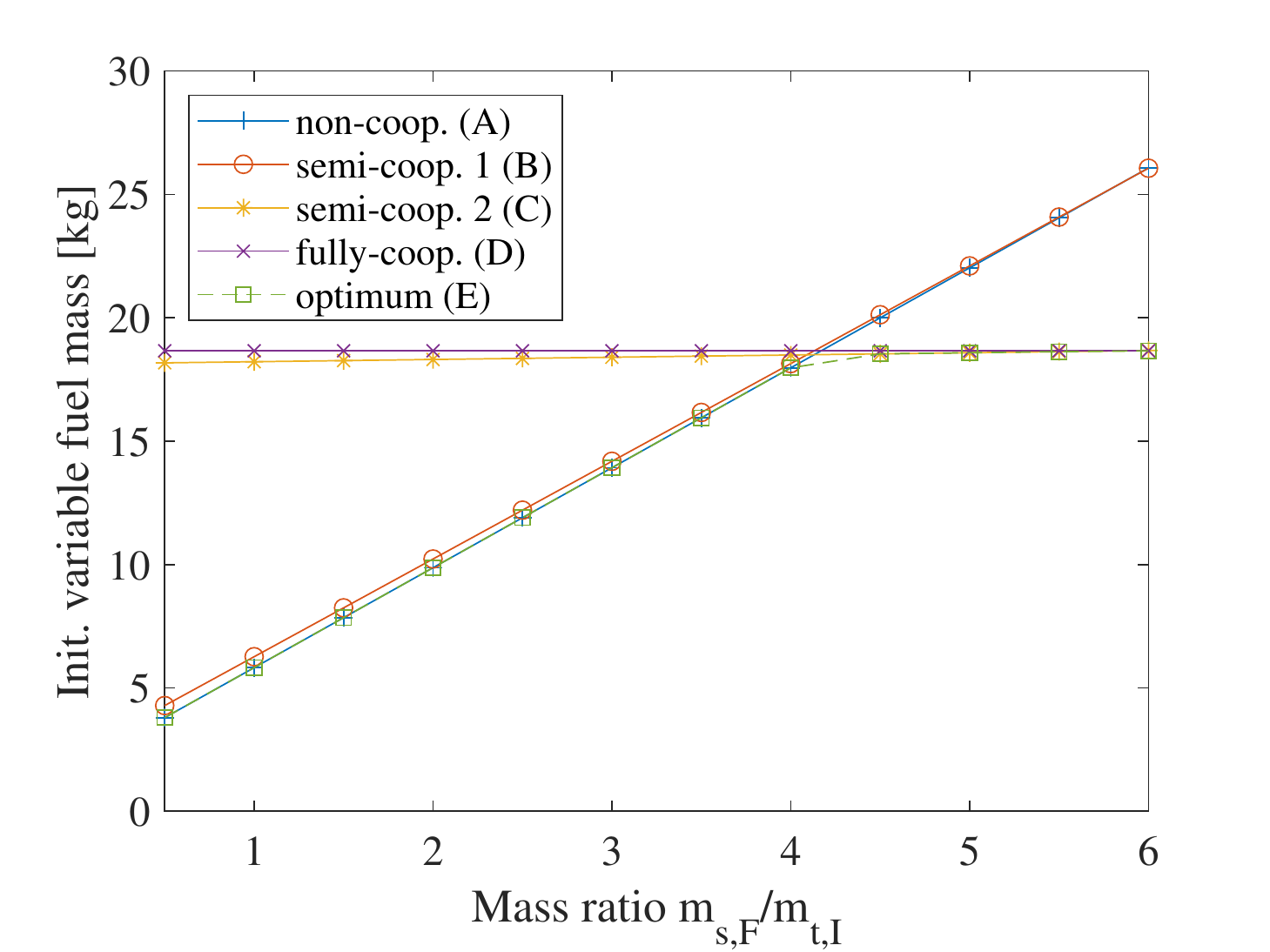}
\caption{$n=12$.}
\label{fig:case_fuelmass_abs_12}
\end{subfigure}
\begin{subfigure}{.49\textwidth}
\centering
\includegraphics[width=.99\textwidth]{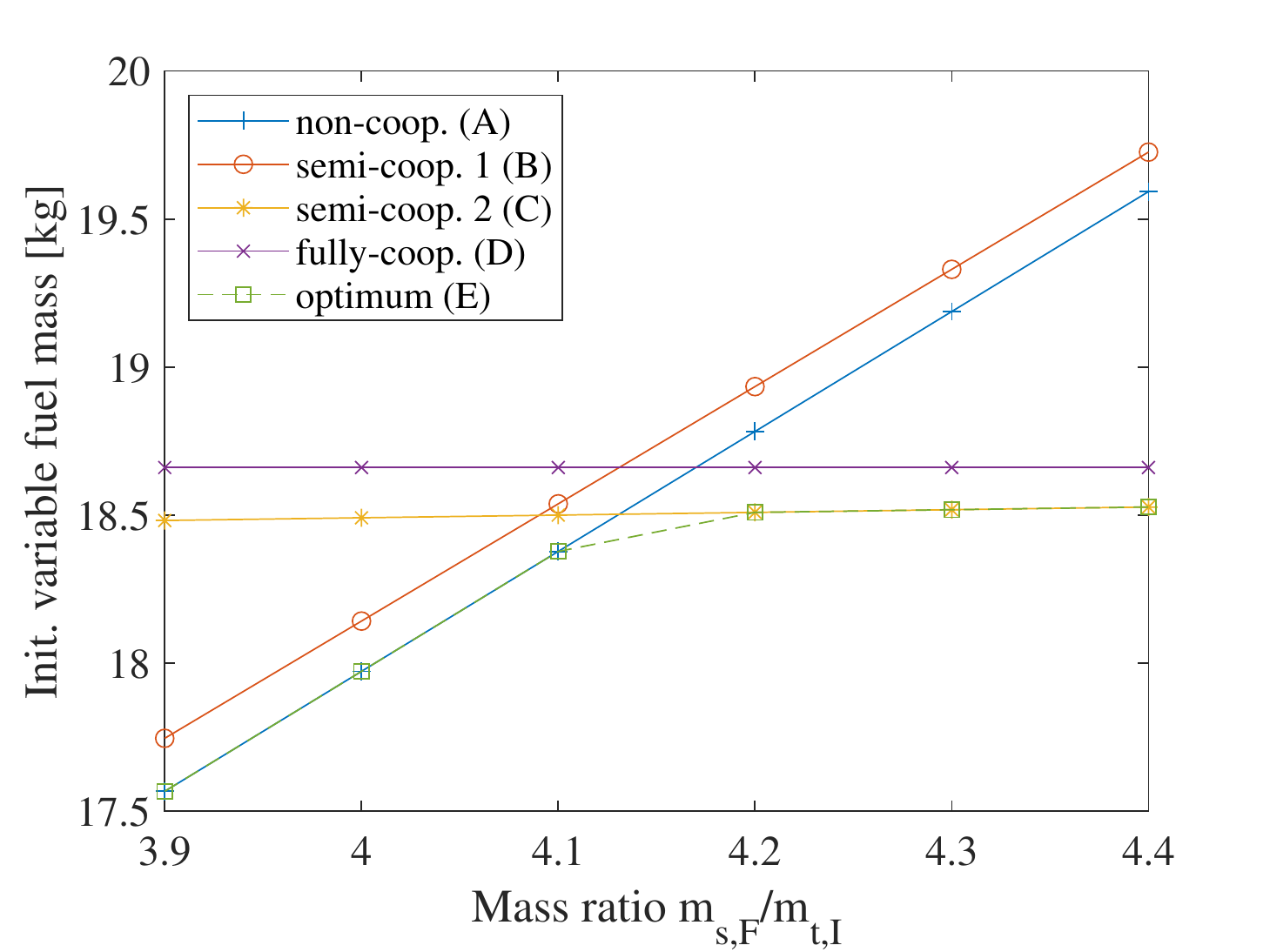}
\caption{$n=12$ around $\alpha_\mathrm{A-D}$.}
\label{fig:case_fuelmass_abs_12_zoom}
\end{subfigure}
\caption{Initial variable fuel mass of each architecture corresponding to the mass ratio $m_\mathrm{s,F}/m_\mathrm{t,I}$. The points where the lines of Architecture A and any other architecture intersect are corresponding to the critical mass ratios $\alpha$.}
\label{fig:case_fuelmass}
\end{figure}

\begin{figure}[hbt!]
\centering
\includegraphics[width=.5\textwidth]{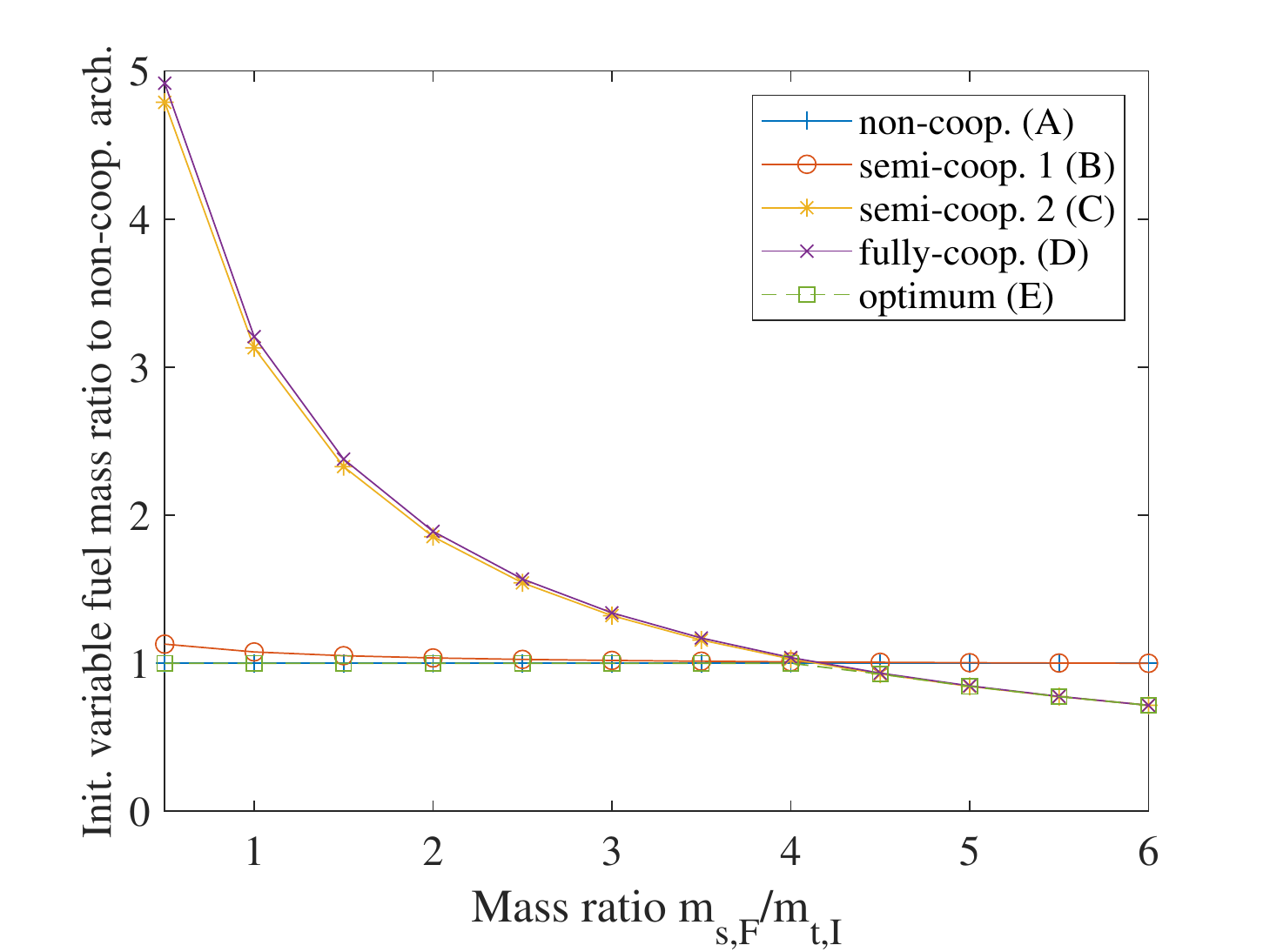}
\caption{Ratio of initial variable fuel mass of each architecture to Architecture A corresponding to the mass ratio $m_\mathrm{s,F}/m_\mathrm{t,I}$.}
\label{fig:case_fuelmass_ratio_12}
\end{figure}

Fig.~\ref{fig:case_fuelmass} compares the initial variable fuel mass in each mission architecture for given mass ratios $m_\mathrm{s,F}/m_\mathrm{t,I}$ for $n=$6, 9, and 12. Note that $n=j \:(j = 1 , 2, \cdots, 12) $ means Targets 1 to $j$ in Fig.~\ref{fig:target_position} are serviced.
The points where the lines of Architecture A and any other architecture intersect are corresponding to the critical mass ratios $\alpha$.
Fig.~\ref{fig:case_fuelmass_abs_12_zoom} shows the initial variable fuel mass of each Architecture for $n = 12$ (Fig.~\ref{fig:case_fuelmass_abs_12}) around the critical mass ratio between Architectures A and D, $\alpha_\mathrm{A-D}$ (see more discussion on $\alpha_\mathrm{A-D}$ in Sec.~\ref{subsubsec:case_crit_mass_ratio}).
Since the servicer does not conduct any orbital transfer in Architecture D, the required variable fuel mass of Architecture D does not depend on the mass of the servicer.
As discussed in Sec.~\ref{sec:analytical_comparison}, when the mass ratio is small (i.e., a light-weight servicer and heavy-weight targets), non-cooperative architecture (Architecture A) requires less fuel mass than the other cooperative architectures. Architecture D with fully cooperative targets requires the most fuel when the mass ratio is small; however, this architecture becomes the most beneficial when the mass ratio is large enough (i.e., a heavy-weight servicer and light-weight targets).
Fig.~\ref{fig:case_fuelmass_ratio_12} compares the initial variable fuel mass as a ratio to that of the non-cooperative Architecture A. As can be seen from this plot, cooperative architectures may require about five times larger amounts of variable fuel compared to the conventional Architecture A when the mass ratio is small while they can also reduce it by 30 \% (or possibly more) when the mass ratio is large.

Furthermore,
% Fig.~\ref{fig:case_fuelmass_abs_12_zoom} shows the initial variable fuel mass of each Architecture for $n = 12$ around the critical mass ratio $\alpha_\mathrm{A-D}$. 
while Architecture E shows the optimal solution with the minimum fuel, this solution can be achieved by either of Architectures A-D in many conditions.
In fact, except for the cases with $n=6$ and $m_\mathrm{s,F}/m_\mathrm{t,I} =$3, 3.5, 4, or 4.5, the optimum architectures turned out to be identical to either Architecture A, C, or D in this case study. In addition, the observed optimum did not decrease the required fuel by a significant amount (see Fig.~\ref{fig:case_fuelmass_abs_6}).
This result suggests that comparing Architectures A, C, and D are practically enough, especially when the mass ratio $m_\mathrm{s,F}/m_\mathrm{t,I}$ is small or large enough compared to the critical mass ratio $\alpha_\mathrm{A-D}$. In other words, mission architects may not have to optimize their mission architecture in these cases.

% \begin{figure}[hbt!]
% \begin{subfigure}{.49\textwidth}
% \centering
% \epsfig {file = figs/fuel_mass_9.eps, width=.99\textwidth}
% \caption{Abs.}
% \label{fig:case_fuelmass_abs_9}
% \end{subfigure}
% \begin{subfigure}{.49\textwidth}
% \centering
% \epsfig {file = figs/fuel_mass_ratio_9.eps, width=.99\textwidth}
% \caption{Ratio.}
% \label{fig:case_fuelmass_rat_9}
% \end{subfigure}
% \caption{Results.}
% \label{fig:case_fuelmass_9}
% \end{figure}

\begin{figure}[hbt!]
\begin{subfigure}{.33\textwidth}
\centering
\includegraphics[width=.99\textwidth]{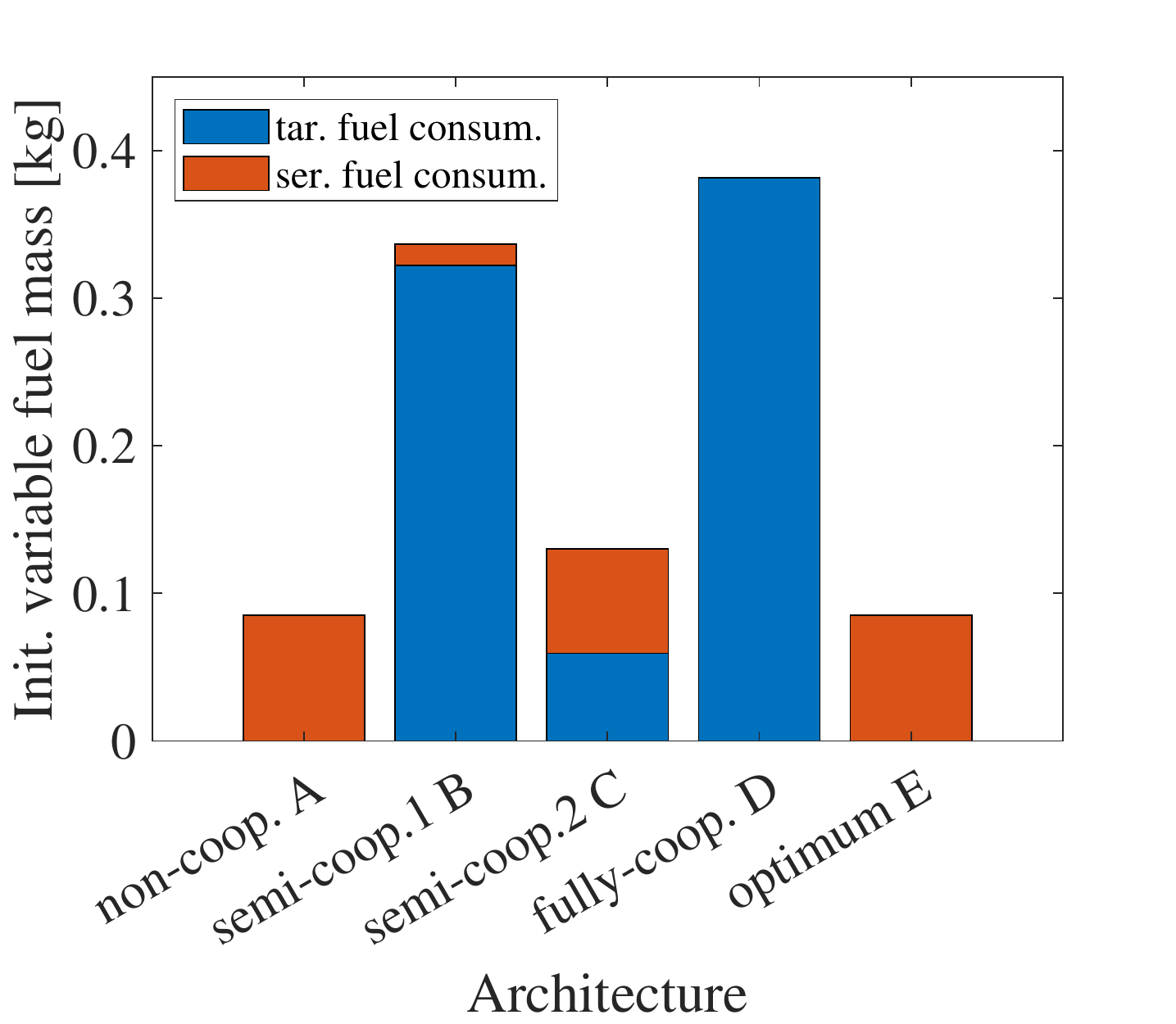}
\caption{$m_\mathrm{s,F}/m_\mathrm{t,I} =0.5$.}
\label{fig:case_fuelmass_bar_0.5_6}
\end{subfigure}
\begin{subfigure}{.33\textwidth}
\centering
\includegraphics[width=.99\textwidth]{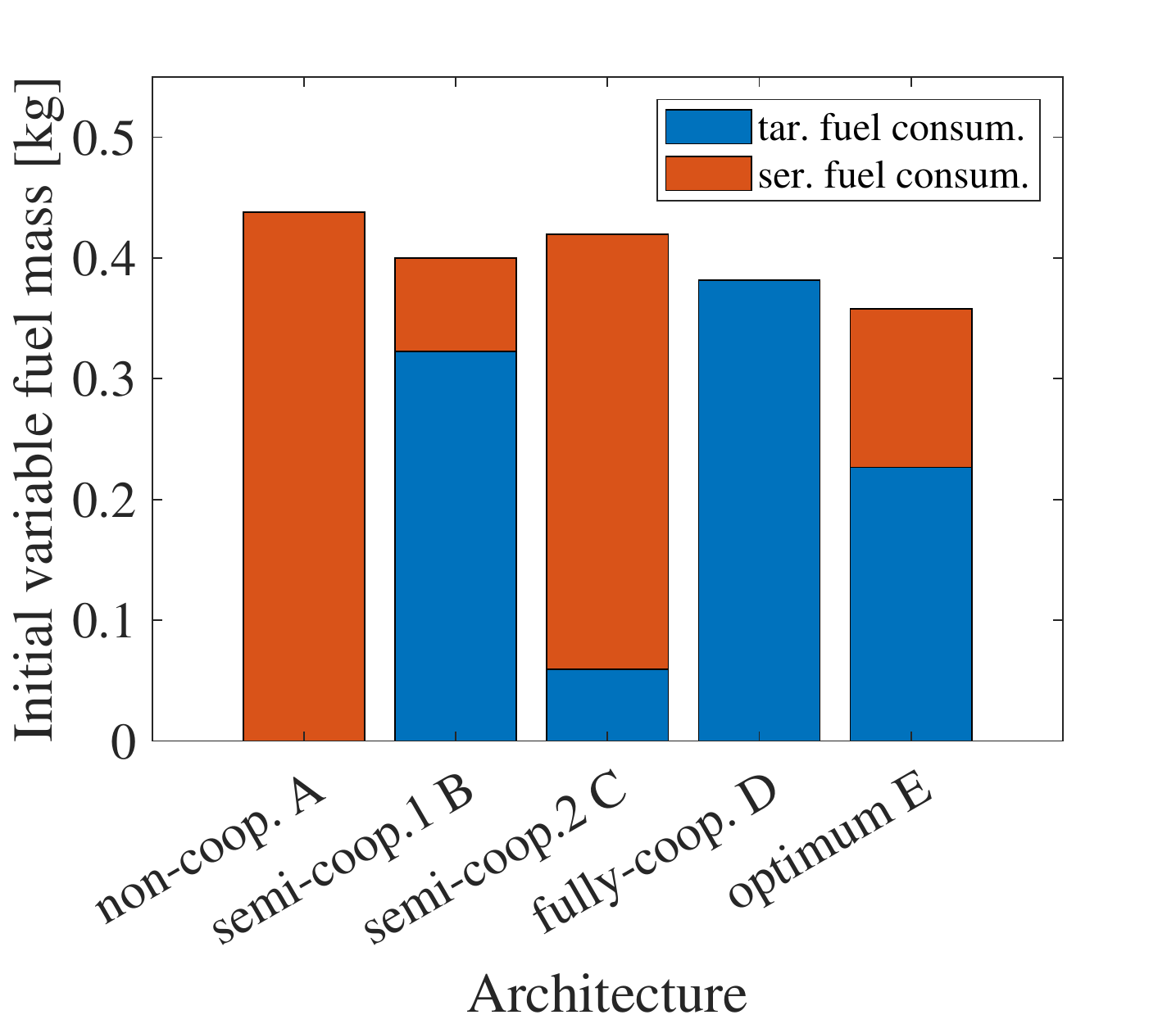}
\caption{$m_\mathrm{s,F}/m_\mathrm{t,I} =4.0$.}
\label{fig:case_fuelmass_bar_4.0_6}
\end{subfigure}
\begin{subfigure}{.33\textwidth}
\centering
\includegraphics[width=.99\textwidth]{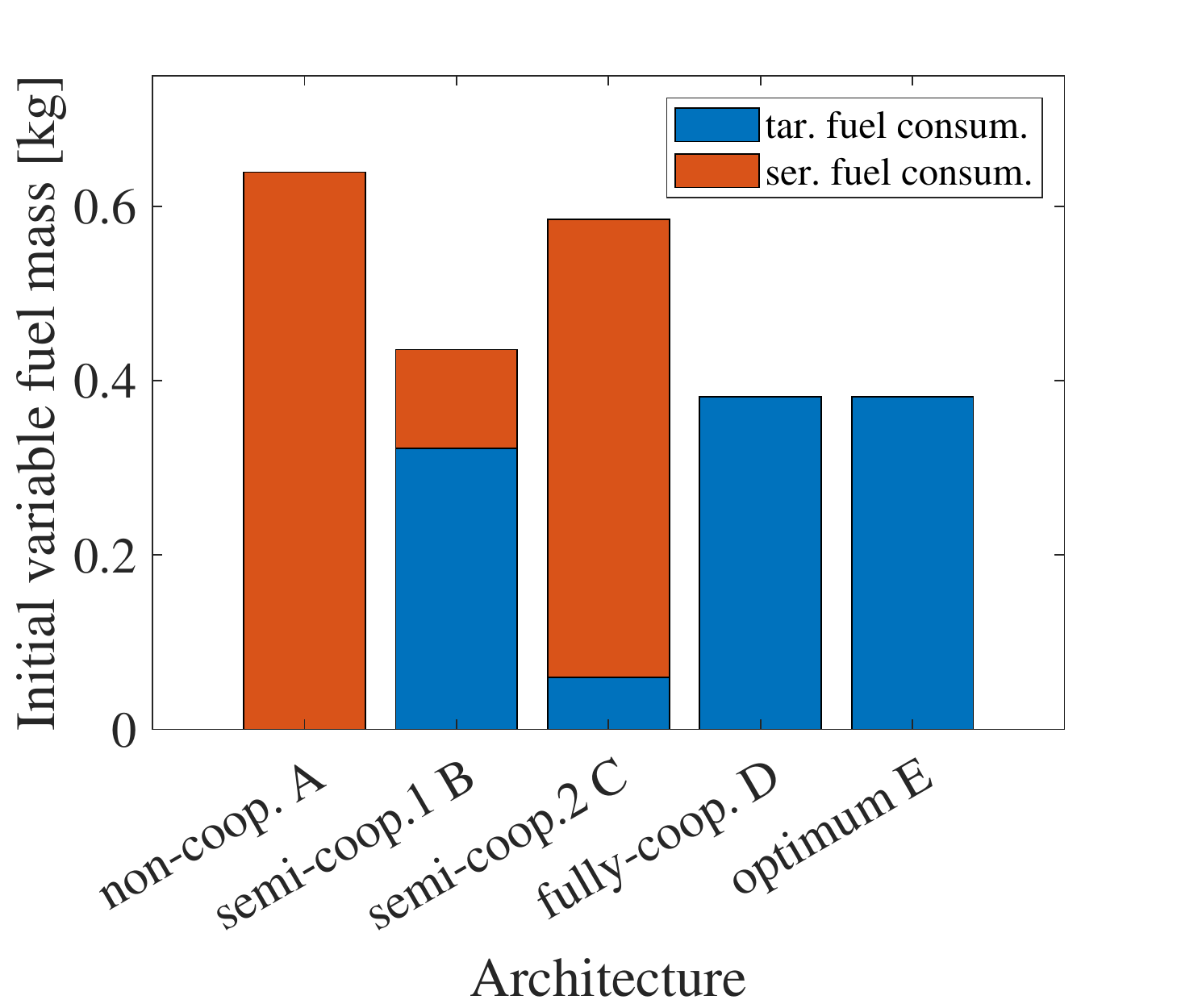}
\caption{$m_\mathrm{s,F}/m_\mathrm{t,I} =6.0$.}
\label{fig:case_fuelmass_bar_6.0_6}
\end{subfigure}
\caption{Initial variable mass breakdown ($n = 6$). Consumed fuel by the targets and the servicer are denoted by "tar. fuel consum." and "ser. fuel consum.," respectively.}
\label{fig:case_fuelmass_bar_n_6}
\end{figure}

\begin{figure}[hbt!]
\begin{subfigure}{.33\textwidth}
\centering
\includegraphics[width=.99\textwidth]{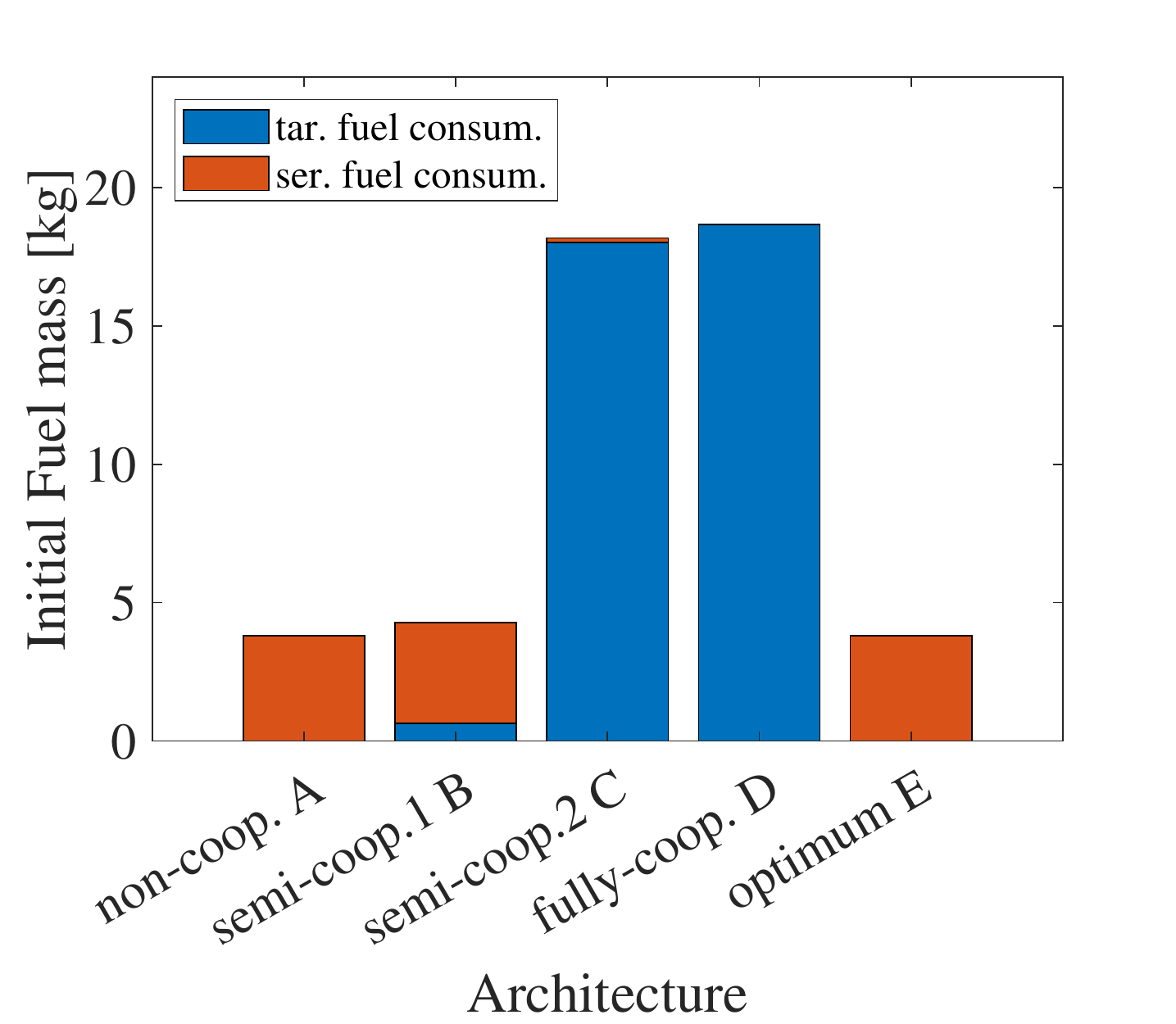}
\caption{$m_\mathrm{s,F}/m_\mathrm{t,I} =0.5$.}
\label{fig:case_fuelmass_bar_0.5}
\end{subfigure}
\begin{subfigure}{.33\textwidth}
\centering
\includegraphics[width=.99\textwidth]{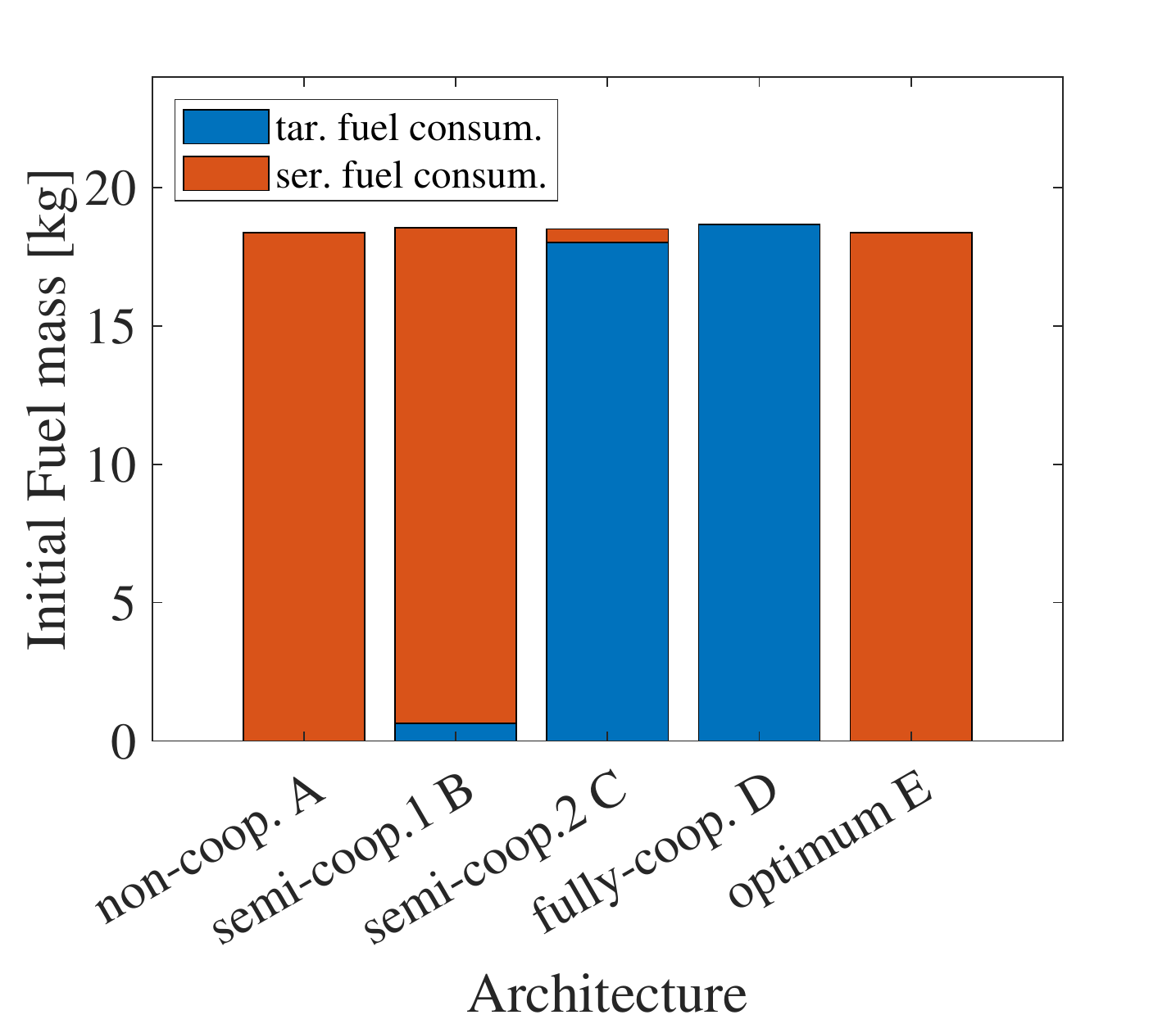}
\caption{$m_\mathrm{s,F}/m_\mathrm{t,I} =4.1$.}
\label{fig:case_fuelmass_bar_4.1}
\end{subfigure}
\begin{subfigure}{.33\textwidth}
\centering
\includegraphics[width=.99\textwidth]{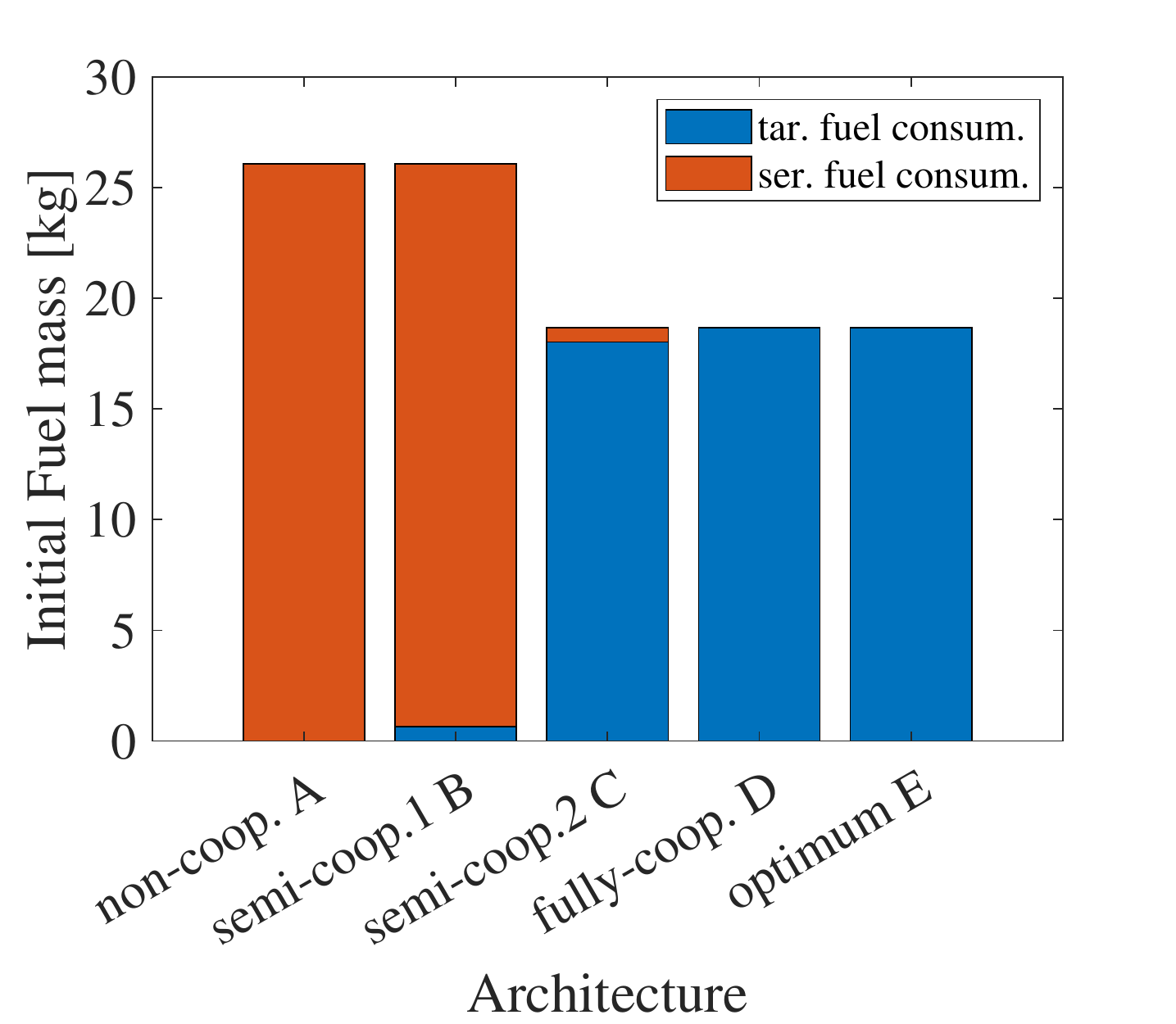}
\caption{$m_\mathrm{s,F}/m_\mathrm{t,I} =6.0$.}
\label{fig:case_fuelmass_bar_6}
\end{subfigure}
\caption{Initial variable mass breakdown ($n = 12$). Consumed fuel by the targets and the servicer are denoted by "tar. fuel consum." and "ser. fuel consum.," respectively.}
\label{fig:case_fuelmass_bar}
\end{figure}

Bar charts Figs.~\ref{fig:case_fuelmass_bar_n_6} and \ref{fig:case_fuelmass_bar} show the breakdown of the initial variable fuel mass of each Architecture with the different mass ratios for the cases of $n = 6$ and 12.
As can be seen from this figure, there's a trade-off between the consumed fuel of the servicer and the targets. This can be particularly seen when focusing on the optimal case in Architecture E. When the servicer final mass is small (Figs.~\ref{fig:case_fuelmass_bar_0.5_6} and \ref{fig:case_fuelmass_bar_0.5}), it is suggested that the servicer actively conducts orbital transfers as in Architecture A since it requires less fuel due to the mass difference between the servicer and the targets.
On the other hand, the targets' active transfers are recommended when the servicer is heavy (Figs.~\ref{fig:case_fuelmass_bar_6.0_6} and \ref{fig:case_fuelmass_bar_6}).
The difference in mass breakdown between Architectures B and C shown in each bar chart of Figs.~\ref{fig:case_fuelmass_bar_n_6} and \ref{fig:case_fuelmass_bar} is due to the gap of the amount of $\Delta v$ required by inclination change and co-planner phasing maneuver.
The relationship between Architectures B and C changes between these two different values of $n$. For instance, Architecture B requires more fuel than C in Fig.~\ref{fig:case_fuelmass_bar_0.5_6} whereas Architecture C requires more in Fig.~\ref{fig:case_fuelmass_bar_0.5}.
This is because the required $\Delta v_\mathrm{inc}$ between the cases $n=6$ and 12 is significantly different due to the very large inclination difference for $n=12$.
This change in the relationship between these two Architectures can also be observed in Fig.~\ref{fig:case_fuelmass} as the slope of each line.

% Comparison between Figs.~\ref{fig:case_fuelmass_bar_n_6} and \ref{fig:case_fuelmass_bar} 
% indicates that the relationship between Architectures B and C changes as the number of the targets changes.
% In the case of $n = 6$, the variable mass fuel breakdown of Architecture B is 

% Figs.~\ref{fig:case_fuelmass_abs_6} and \ref{fig:case_fuelmass_abs_9} 

\subsubsection{Critical Mass Ratio}
\label{subsubsec:case_crit_mass_ratio}
\begin{figure}[hbt!]
\centering
\includegraphics[width=.7\textwidth]{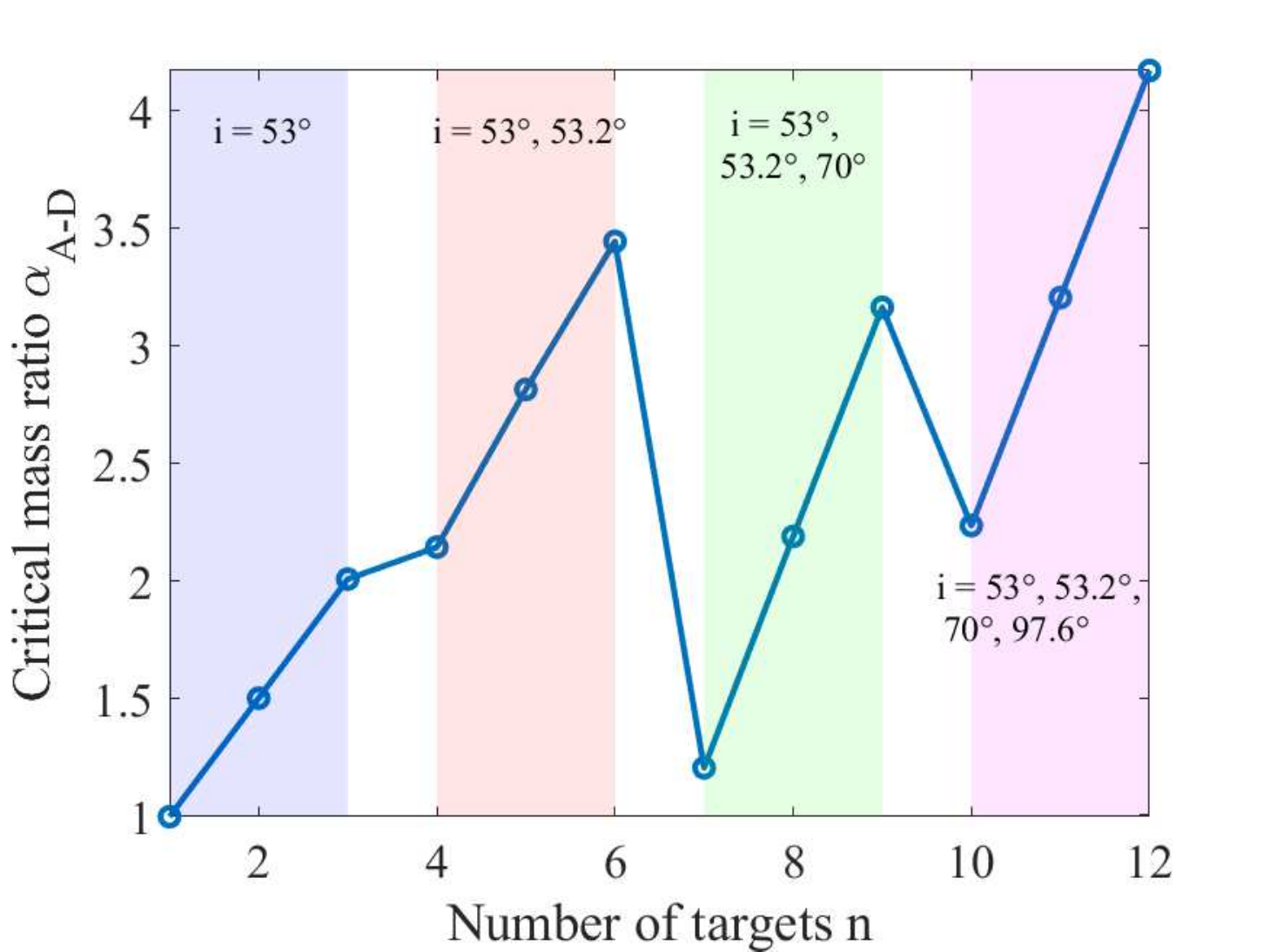}
\caption{Change in the critical mass ratio between Architectures A and D in response to various numbers of targets $n$. Targets 1 to $n$ are serviced when $n \neq 1$.}
\label{fig:mass_critical_case}
\end{figure}

We next examine further the critical mass ratio, particularly the one that differentiates Architectures A and D. For Architecture D with fully cooperative targets, the servicer does not conduct any orbital transfer, i.e., $\Delta v_\mathrm{s,c} = 0$. 
Therefore, from 
% Eq.~\eqref{eq:mass_critical},
Inequality \eqref{eq:mass_ineq2},
the critical mass ratio that differentiates Architectures A and D, $\alpha_\mathrm{A-D}$, can be calculated as
\begin{equation}
\label{eq:mass_critical_A-D}
    \alpha_\mathrm{A-D} = \frac{ \sum_{j = 1}^{n} \left\{ \left( \frac{m_\mathrm{req}}{m_\mathrm{t,I}} \exp\left(\frac{\sum_{k = 1}^{j}\Delta v_\mathrm{s,n} (k)}{I_\mathrm{sp, s} g_0}\right) \right) -   \left[\left(1 +\frac{m_\mathrm{req}}{m_\mathrm{t,I}}  \right) \exp\left(\frac{\Delta v_\mathrm{t, out} (j)}{I_\mathrm{sp, t} g_0}\right) - \exp\left(\frac{-\Delta v_\mathrm{t, in} (j)}{I_\mathrm{sp, t} g_0}\right)\right] \right\}  }{1 - \exp\left(\frac{\sum_{j = 1}^{n+1}\Delta v_\mathrm{s,n} (j)}{I_\mathrm{sp, s} g_0}\right)}.
\end{equation}
For each number of targets, $\alpha_\mathrm{A-D}$ is calculated using Eq.~\eqref{eq:mass_critical_A-D}, and the results are shown in Fig.~\ref{fig:mass_critical_case}. 
% Note that $n=j \:(j = 1 , 2, \cdots, 12) $ means Targets 1 to $j$ in Fig.~\ref{fig:target_position} are serviced.
When $n=1$, $\alpha_\mathrm{A-D}$ becomes 1 suggesting whichever lighter between the servicer and Target 1 should move to save some fuel.

\begin{figure}[hbt!]
\center
\includegraphics[width=.85\textwidth]{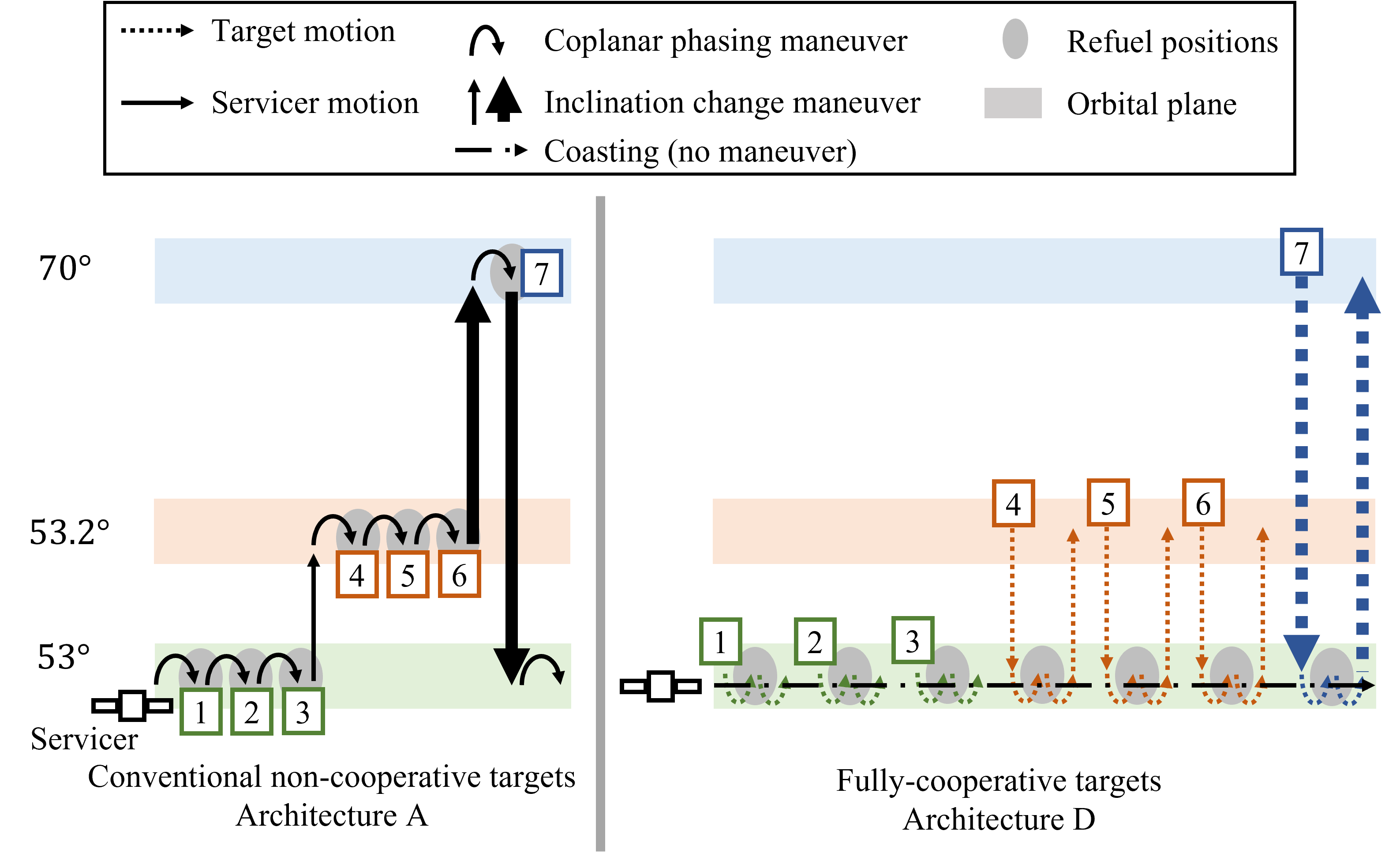}
\caption{Comparison of mission architectures ($n=7$).}
\label{fig:architecture_n_7}
\end{figure}

As can be seen in Fig.~\ref{fig:mass_critical_case}, $\alpha_\mathrm{A-D}$ linearly increases when the number of targets increases in the same orbital planes (see each colored area in this figure). When a new target in a different orbital plane is added, $\alpha_\mathrm{A-D}$ can decrease. For instance, there is a large drop between $n=$ 6 and 7, and $\alpha_\mathrm{A-D}$ becomes close to 1.0 when $n=7$.
This is due to the large $\Delta v_\mathrm{inc}$ required for the refuel to Target 7. Fig.~\ref{fig:architecture_n_7} compares Architecture A and D when $n=7$.
Since the required $\Delta v_\mathrm{inc}$ before and after the 7th refuel is significantly larger in both cases than that for other coplanar phasing maneuvers and inclination changes between $i= 53^{\circ}$ and $53.2^{\circ}$ (see bold arrows in Fig.~\ref{fig:architecture_n_7}), these  orbital transfers become the main mass driver. Since the difference between Architectures A and D regarding these transfers is rather small (i.e., $i=53.2^{\circ}$ to $70^{\circ}$ to $53^{\circ}$ for Architecture A, and $i=70^{\circ}$ to $53^{\circ}$ to $70^{\circ}$ for Architecture D), the critical mass ratio becomes close to 1.0.
The behavior of the critical mass ratio is further discussed in Sec.~\ref{subsec:case_sensitivity}.

% \begin{table}[hbt!]
% \center
% \caption{Critical mass ratio}
% \label{tab:mass_critical_case}
% \begin{tabular}{rccccccccccccc}
% \thickhline
% \# of targets& $n$ & 1 & 2 & 3 & 4 & 5& 6 & 7& 8 & 9 & 10 & 11& 12\\\hline
% critical mass ratio &$\alpha_\mathrm{A-D}$ & 1.00 & 1.50 & 2.01 & 2.14 & 2.82 & 3.44 & 1.21 & 2.19 & 3.16 & 2.24 & 3.21 & 4.17\\\hline
% \end{tabular}
% \end{table}

% 0.999993064569712	1.50344636072649	2.00922720444099	2.14493086261843	2.81524957874011	3.44399998572558	1.20747315445176	2.19063824447137	3.16372975149034	2.23650667267260	3.20529373044255	4.16995129248172

\subsection{Discussion and Sensitivity Analysis}
\label{subsec:case_sensitivity}
As previously mentioned in Sec.~\ref{subsec:case_init_mass}, Fig.~\ref{fig:mass_critical_case} suggests that the critical mass ratio generally increases as the number of targets increases in the same orbital plane; however, it does not necessarily increase when the newly added targets are in a different orbital plane. To study this behavior, in this subsection, the sensitivity of inclination angles on
the critical mass ratio is first examined by employing target sets where targets are distributed differently from Fig.~\ref{fig:target_position}. Table \ref{tab:target_sets} lists target sets considered in this analysis and their inclination angles. Note that each target set has the same set of the argument of latitude $u$ values as the original Target Set a (see Fig.~\ref{fig:target_position}).
These target sets are selected to explore the effects of multi-plane refueling compared to coplanar refueling (Target Set b and the other sets) and the large inclination change both near the initial orbit of the servicer (Target Sets c and d)  and far from that (Target Sets a and e).

\begin{table}[hbt!]
\center
\caption{Target sets and their inclination angles [$^{\circ}$].}
\label{tab:target_sets}
\begin{tabular}{lcccc}
\thickhline
Target \#                  & 1-3 & 4-6  & 7-9  & 10-12 \\\hline
a (original) & 53  & 53.2 & 70   & 97.6  \\
                             b            & 53  & 53   & 53   & 53    \\
                             c            & 53  & 53.2 & 53.2 & 53.2  \\
                             d            & 53  & 70   & 70   & 70    \\
                             e            & 53  & 53.2 & 70   & 75   \\\hline
\end{tabular}
\end{table}

\begin{figure}[hbt!]
\centering
\includegraphics[width=.7\textwidth]{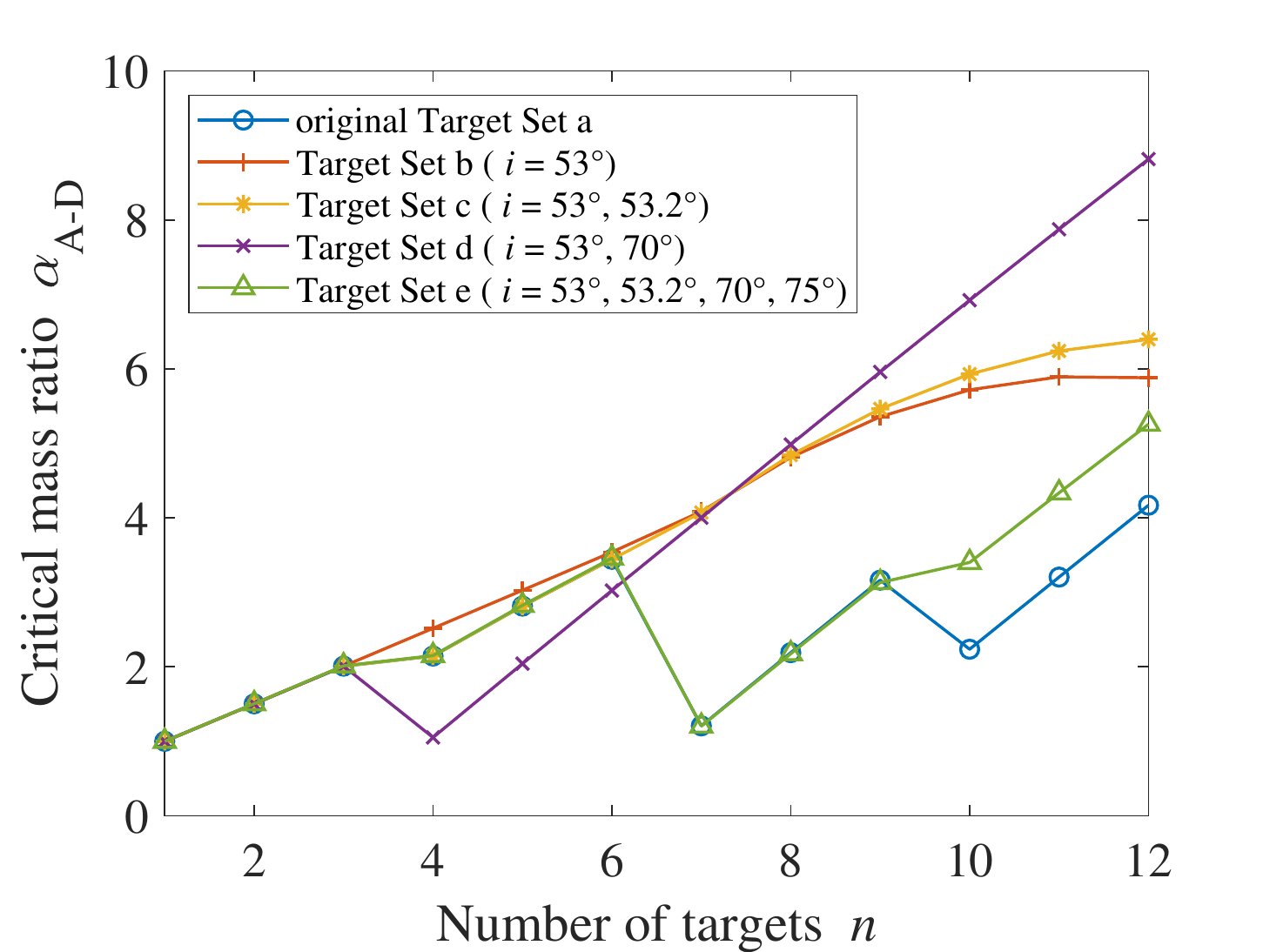}
\caption{Effect of different target sets on the critical mass ratio.}
\label{fig:mass_critical_case_diff_plane}
\end{figure}

Fig.~\ref{fig:mass_critical_case_diff_plane} illustrates the behavior of the critical mass ratio in response to the change in $n$ with different target sets listed in Table \ref{tab:target_sets}.
Since all targets in Target Set b are coplanar, the trend of this red line does not change even if all targets are in either $i=70^{\circ}$, $97.6^{\circ}$, or any other plane.
As expected, large drops of $\alpha_\mathrm{A-D}$ occur when a new Target $j$ is added to a new plane that is far from the plane where the Target $(j-1)$  is initially located. On the other hand, this result suggests that when the newly-added plane of Target $j$ is close to the plane of the Target $(j-1)$, the effect of adding a new plane becomes less significant. These behaviors can be observed by comparing the lines of Target Sets b, c, and d around $n=3$ and 4, and the lines of Target Sets a and e around $n=9$ and 10.
A comparison of the lines of Target Sets c and d also suggests that the increment of the critical mass ratio becomes larger when the difference of inclination between the initial plane and the newly added plane is larger.

Furthermore, for the coplanar case (Target Set b), $\alpha_\mathrm{A-D}$ for $n=11$ and 12 are almost the same value. This is because the newly added Target 12 is located close to the original position of the servicer indicating that this orbit is becoming crowded with targets and adding one more target affects the entire refueling architectures less significantly.
% the newly required $\Delta v$s for Architecture D, $\Delta v_\mathrm{t, in}(12)$ and $\Delta v_\mathrm{t, out}(12)$, are relatively small as well as the change in required $\Delta v_\mathrm{s, n}$ in Architecture A. Therefore, the change in required $\Delta v$ between $n=11$ and $12$ is small and similar between Architecture A and D, which makes the critical mass ratio for these cases almost the same.
The line of Target Set c shows a similar trend: the slope of the Target Set c line becomes less steep around $n=12$. This is because target set c is almost coplanar whereas the orbital planes of Target Sets a, d, and e are more scattered.

% \begin{figure}[hbt!]
% \centering
% \includegraphics[width=.99\textwidth]{figs/architecture_n_12_coplanar.png}
% \caption{Comparison of mission architectures with coplanar refueling ($n=11$ and $12$).}
% \label{fig:architecture_n_12_coplanar}
% \end{figure}

\begin{figure}[hbt!]
\begin{subfigure}{.49\textwidth}
\centering
\includegraphics[width=.99\textwidth]{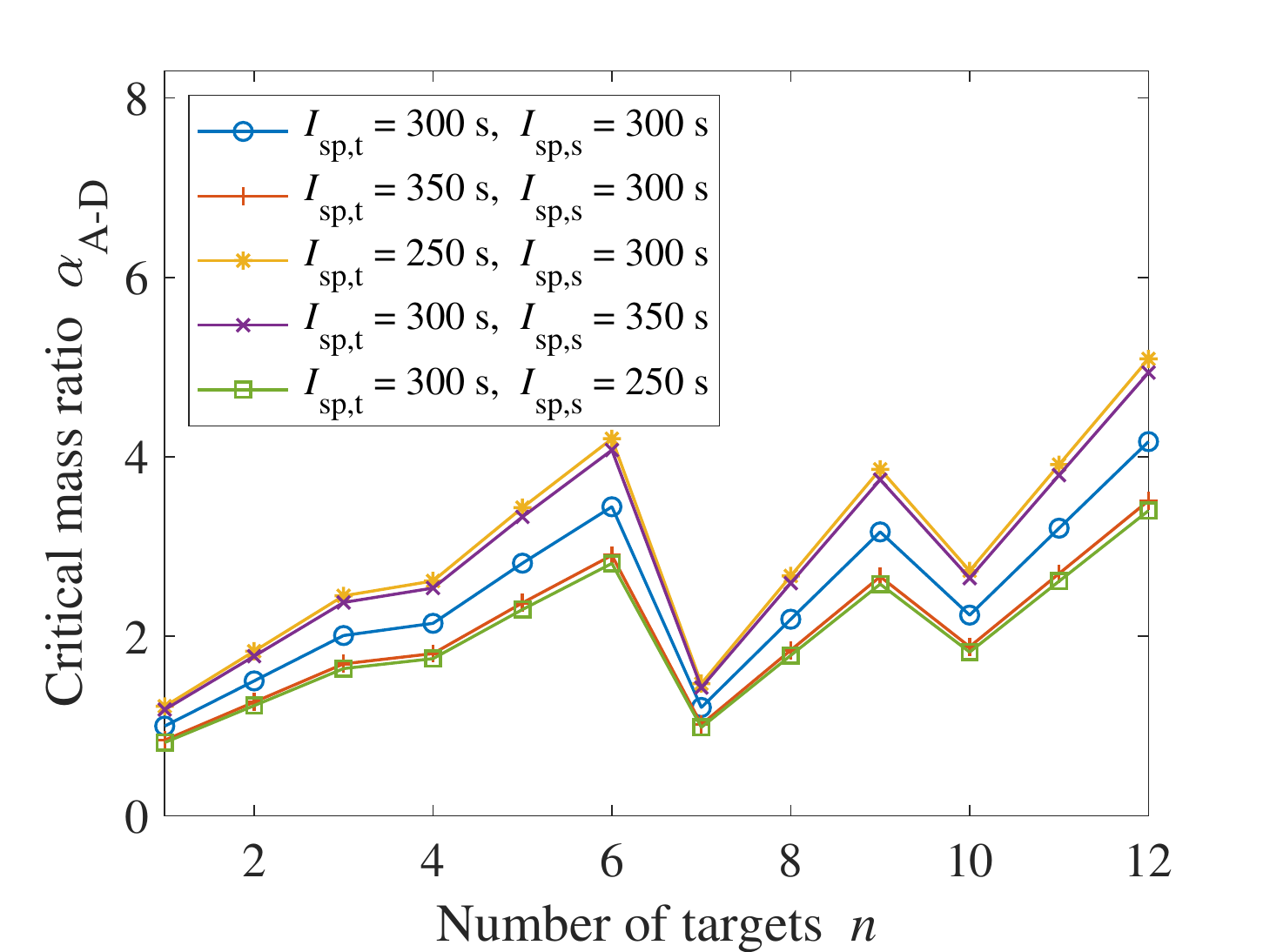}
\caption{Target Set a (multi-plane refueling).}
\label{fig:mass_critical_case_diff_isp_a}
\end{subfigure}
\begin{subfigure}{.49\textwidth}
\centering
\includegraphics[width=.99\textwidth]{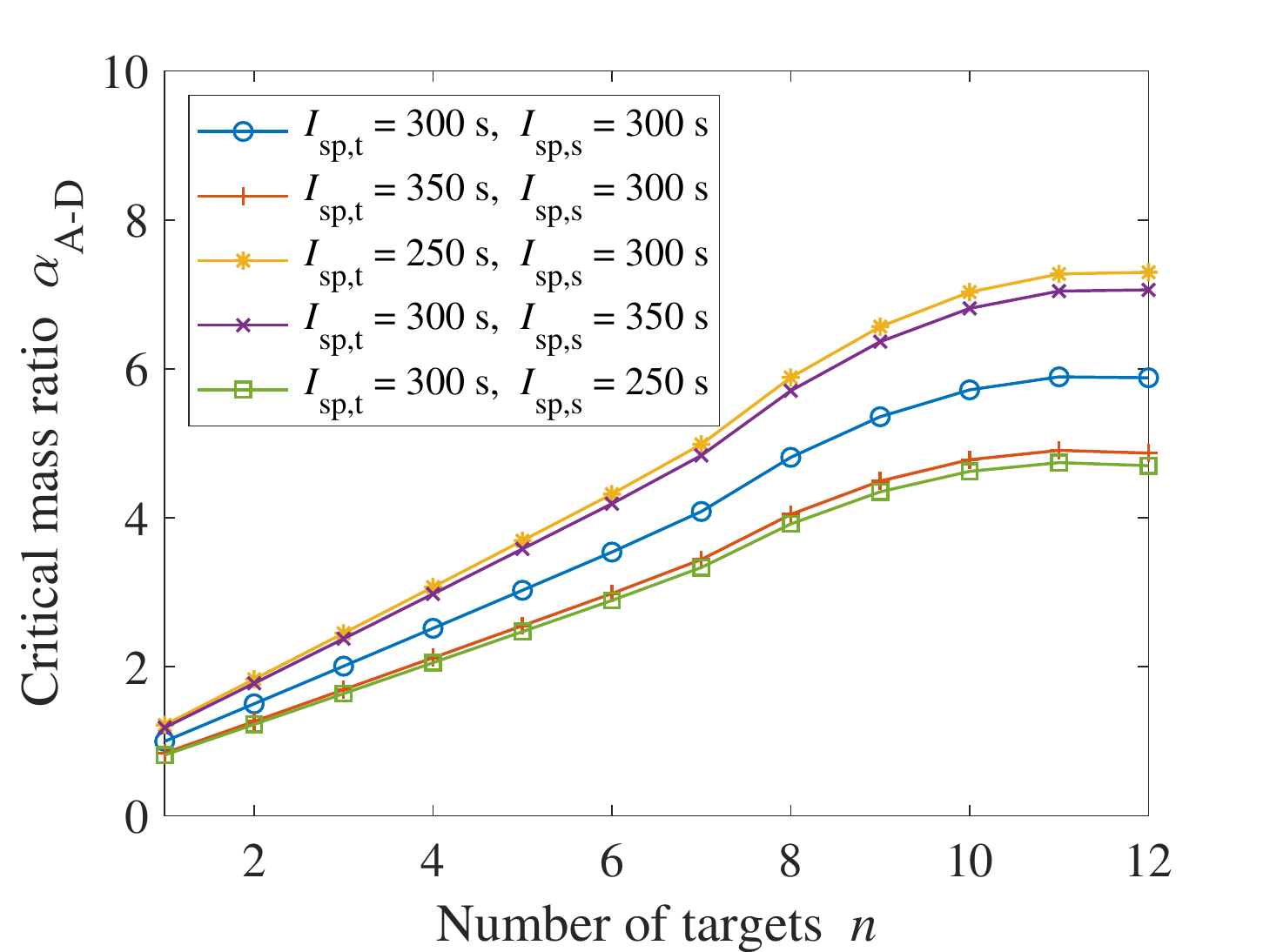}
\caption{Target Set b (coplanar refueling).}
\label{fig:mass_critical_case_diff_isp_b}
\end{subfigure}
\caption{Effect of different $I_\mathrm{sp}$ sets on the critical mass ratio.}
\label{fig:mass_critical_case_diff_isp}
\end{figure}

% finally if every 12 target is in different orbital planes
We also compare the effect of the specific impulse of the target ($I_\mathrm{sp, t}$) and the servicer ($I_\mathrm{sp, s}$). Eq.~\eqref{eq:mass_critical_A-D} implies a smaller $I_\mathrm{sp, s}$ or a larger $I_\mathrm{sp, s}$ leads to a larger $\alpha_\mathrm{A-D}$, and conversely, a larger $I_\mathrm{sp, s}$ or a smaller $I_\mathrm{sp, s}$ leads to a smaller $\alpha_\mathrm{A-D}$.
Fig.~\ref{fig:mass_critical_case_diff_isp} shows the results of employing different values for both $I_\mathrm{sp, t}$ and $I_\mathrm{sp, s}$ for Target Set a (Fig.~\ref{fig:mass_critical_case_diff_isp_a}) and Target Set b (Fig.~\ref{fig:mass_critical_case_diff_isp_b}).
As can be seen from this figure, the sensitivity of each specific impulse turned out to be about the same level between both target sets.
Moreover, as expected from Eq.~\eqref{eq:mass_critical_A-D}, increasing $I_\mathrm{sp,t}$ or decreasing $I_\mathrm{sp,s}$ reduces the critical mass ratio implying the fully-cooperative case can be favored with smaller servicer mass.
On the other hand, smaller $I_\mathrm{sp,t}$ or larger $I_\mathrm{sp,s}$ makes the conventional non-cooperative architecture more favored.
The sensitivities of $I_\mathrm{sp, t}$ and $I_\mathrm{sp, s}$ also turn out to be the same according to this figure.
The investigation of the effect of employing electric propulsion is left for future work since the continuous orbital transfer significantly complicates the operation of this on-orbit refueling.

% \begin{figure}[hbt!]
% \centering
% \epsfig {file = figs/alpha_diff_isp.eps, width=.7\textwidth}
% \caption{Effect of different $I_\mathrm{sp}$ sets on the critical mass ratio.}
% \label{fig:mass_critical_case_diff_isp}
% \end{figure}

% Eq.~\eqref{eq:mass_critical_A-D}, suggest that the ratio of required refuel mass to the initial target mass, $m_\mathrm{req}/m_\mathrm{t,I}$, can be another parameter that affects the critical mass ratio.
% However, as Fig.~\ref{fig:mass_critical_case_fuel_ratio} indicates this ratio does not affect the critical mass ratio as much as the number of targets, and therefore, mission architects can roughly assume this value from the target spacecraft design.
% Note that this ratio can be larger than 1.0, since $m_\mathrm{t,I}$ denotes the mass of a target at the beginning of the whole refueling campaign but not the mass at the initial orbital insertion.

% % \begin{figure}[hbt!]
% % \centering
% % \epsfig {file = figs/alpha_fuel_ratio.eps, width=.7\textwidth}
% % \caption{Effect of refuel mass ratio $m_\mathrm{req}/m_\mathrm{t,I}$ on the critical mass ratio.}
% % \label{fig:mass_critical_case_fuel_ratio}
% % \end{figure}

% \begin{figure}[hbt!]
% \begin{subfigure}{.49\textwidth}
% \centering
% \epsfig {file = figs/alpha_fuel_ratio_a.eps, width=.99\textwidth}
% \caption{Target set a (multi-plane refueling).}
% % \label{fig:case_fuelmass_abs_6}
% \end{subfigure}
% \begin{subfigure}{.49\textwidth}
% \centering
% \epsfig {file = figs/alpha_fuel_ratio_b.eps, width=.99\textwidth}
% \caption{Target set b (coplanar refueling).}
% % \label{fig:case_fuelmass_abs_9}
% \end{subfigure}
% \caption{Effect of refuel mass ratio $m_\mathrm{req}/m_\mathrm{t,I}$ on the critical mass ratio.}
% \label{fig:mass_critical_case_fuel_ratio}
% \end{figure}

\section{Conclusion}
\label{sec:conclusion}
This paper developed an analytical model for non-cooperative and cooperative refueling architectures, and analytically derived the critical mass ratio that represents the break-even point between the two architectures. We further analyzed how mission parameters change this value through the case study of non-coplanar multi-target on-orbit refueling.
If the ratio between the final servicer mass and the initial target mass,~$m_\mathrm{s,F}/m_\mathrm{t,I}$, is smaller than the critical mass ratio (a light-weight servicer and heavy-weight targets), the conventional non-cooperative architecture requires less fuel than the cooperative architecture. Conversely, if this ratio is larger than the critical mass ratio (a heavy-weight servicer with light-weight targets), employing the cooperative architecture may save some fuel.
The result from the case study suggests this critical mass ratio generally increases as the number of targets increases in the same orbital plane. However, when new targets locate in different orbital planes, the critical mass ratio can also decrease depending on the difference in the inclination.
The optimization of rendezvous inclination and the argument of latitude suggests that the saved fuel compared to non-cooperative architecture or fully-cooperative architecture is rather small. Hence, it may be sufficient to consider either of these architectures, especially for the initial mission design.

The extension of this study in spacecraft with low-thrust and high-specific-impulse propulsion (e.g., solar electric propulsion) is left for future work.
The optimization of the order of spacecraft to be refueled and the study of the effect of refueling architecture on the order is another direction of future work.
Different transfer techniques between each refuel such as solving Lambert's problems should be discussed in the extension of this work, too.

\appendix
\section*{Appendix}
\label{sec:appendix}
\subsection{Analytical Expression of Critical Mass Ratio}
\label{app:crit_mass_ratio}

From Eqs.\eqref{eq:mass_final_coop} and \eqref{eq:mass_final_noncoop}, the critical mass ratio can be derived as follows:
\begin{align}
    \begin{split}
        &{}m_\mathrm{s,F} \exp\left(\frac{\sum_{j = 1}^{n+1}\Delta v_\mathrm{s,c} (j)}{I_\mathrm{sp, s} g_0}\right) + \sum_{j = 1}^{n} \left\{ \left[\left( m_\mathrm{t,I} + m_\mathrm{req} \right) \exp\left(\frac{\Delta v_\mathrm{t, out} (j)}{I_\mathrm{sp, t} g_0}\right) - m_\mathrm{t,I}  \exp\left(\frac{-\Delta v_\mathrm{t, in} (j)}{I_\mathrm{sp, t} g_0}\right)\right] \exp\left(\frac{\sum_{k = 1}^{j}\Delta v_\mathrm{s,c} (k)}{I_\mathrm{sp, s} g_0}\right) \right\}\nonumber\\
        &{}=   m_\mathrm{s,F} \exp\left(\frac{\sum_{j = 1}^{n+1}\Delta v_\mathrm{s,n} (j)}{I_\mathrm{sp, s} g_0}\right) + \sum_{j = 1}^{n} \left(m_\mathrm{req} \exp\left(\frac{\sum_{k = 1}^{j}\Delta v_\mathrm{s,n} (k)}{I_\mathrm{sp, s} g_0}\right) \right)\nonumber
    \end{split}\\
    \begin{split}
    % \label{eq:mass_ineq2}
        \Leftrightarrow &{}\frac{m_\mathrm{s,F}}{m_\mathrm{t,I}} \left(\exp\left(\frac{\sum_{j = 1}^{n+1}\Delta v_\mathrm{s,c} (j)}{I_\mathrm{sp, s} g_0}\right) - \exp\left(\frac{\sum_{j = 1}^{n+1}\Delta v_\mathrm{s,n} (j)}{I_\mathrm{sp, s} g_0}\right)\right)\\
        &{}= \sum_{j = 1}^{n} \left\{ \left( \frac{m_\mathrm{req}}{m_\mathrm{t,I}} \exp\left(\frac{\sum_{k = 1}^{j}\Delta v_\mathrm{s,n} (k)}{I_\mathrm{sp, s} g_0}\right) \right) -   \left[\left(1 +\frac{m_\mathrm{req}}{m_\mathrm{t,I}}  \right) \exp\left(\frac{\Delta v_\mathrm{t, out} (j)}{I_\mathrm{sp, t} g_0}\right) - \exp\left(\frac{-\Delta v_\mathrm{t, in} (j)}{I_\mathrm{sp, t} g_0}\right)\right] \exp\left(\frac{\sum_{k = 1}^{j}\Delta v_\mathrm{s,c} (k)}{I_\mathrm{sp, s} g_0}\right)\right\}\nonumber
    \end{split}\\
    \begin{split}
    % \label{eq:mass_ineq2}
        \Leftrightarrow &{}\frac{m_\mathrm{s,F}}{m_\mathrm{t,I}} =  \frac{ \sum_{j = 1}^{n} \left\{ \left( \frac{m_\mathrm{req}}{m_\mathrm{t,I}} \exp\left(\frac{\sum_{k = 1}^{j}\Delta v_\mathrm{s,n} (k)}{I_\mathrm{sp, s} g_0}\right) \right) -   \left[\left(1 +\frac{m_\mathrm{req}}{m_\mathrm{t,I}}  \right) \exp\left(\frac{\Delta v_\mathrm{t, out} (j)}{I_\mathrm{sp, t} g_0}\right) - \exp\left(\frac{-\Delta v_\mathrm{t, in} (j)}{I_\mathrm{sp, t} g_0}\right)\right] \exp\left(\frac{\sum_{k = 1}^{j}\Delta v_\mathrm{s,c} (k)}{I_\mathrm{sp, s} g_0}\right)\right\}  }{\exp\left(\frac{\sum_{j = 1}^{n+1}\Delta v_\mathrm{s,c} (j)}{I_\mathrm{sp, s} g_0}\right) - \exp\left(\frac{\sum_{j = 1}^{n+1}\Delta v_\mathrm{s,n} (j)}{I_\mathrm{sp, s} g_0}\right)}.
    \end{split}
\end{align}

\bibliography{sample}

\end{document}